\documentclass[12pt,leqno]{article}
\usepackage{amssymb,amsmath,amsfonts,amsthm,color,mathrsfs}
\usepackage{bbm}
\usepackage{dsfont}
\usepackage[marginal]{footmisc}
\usepackage[]{hyperref}
\usepackage{color} 
\usepackage{mathtools}
\usepackage{graphicx}
\usepackage{enumitem}
\theoremstyle{plain}

\newtheorem{theorem}{Theorem}[section]
\newtheorem{lemma}[theorem]{Lemma}
\newtheorem{remark}[theorem]{Remark}
\newtheorem{corollary}[theorem]{Corollary}
\theoremstyle{definition}
\newtheorem{definition}[theorem]{Definition}

\newtheorem{assumption}[theorem]{Assumption}

\definecolor{wco}{rgb}{0.5,0.2,0.3}
\allowdisplaybreaks[4]

\numberwithin{equation}{section} \theoremstyle{remark}

\def\R{\mathbb{R}}
\def\Rd{\mathbb{R}^d}
\def\E{\mathbb{E}}
\def\P{\mathscr{P}}

\def\L{\mathscr{L}}
\def\M{\mathcal{M}}
\def\d{\textup{d}}
\def\D{\mathcal{D}}

\def\tr{\textup{tr}}

\def\<{\langle}
\def\>{\rangle}
\def\Proof.{\noindent{\bf Proof. }}

\def\e{\epsilon}

\def\div{\textup{div}}

\def\loc{\textup{loc}}

\def\newdot{{\kern.8pt\cdot\kern.8pt}}

\def \b{\mathfrak{b} }
\def\log{\operatorname{log}}
\begin{document}
	
	\title{\bf Well-posedness and Quantitative Stability for Degenerate Fokker--Planck Equations with Rough Coefficients}
	
	\author{
		\begin{tabular}{cc}
			{\bf Huaiqian Li$^{a}$} & {\bf Liying Mu$^{a}$\thanks{Corresponding author.}} \\
			\footnotesize huaiqian.li@tju.edu.cn & \footnotesize liying\_mu01@tju.edu.cn
		\end{tabular} \\ \vspace{0.5em}
		\footnotesize{$^{a}$ Center for Applied Mathematics and KL-AAGDM, Tianjin University, Tianjin 300072, P. R. China}
	}
	
	\date{}

	\maketitle
\begin{abstract}
	We establish quantitative stability and well-posedness for degenerate Fokker--Planck equations with Sobolev-regular and possibly fully degenerate diffusion coefficient. For two classes of rough drifts, stability estimates imply uniqueness, and existence is proved under suitable conditions. First, for anisotropic drifts with partial fractional Sobolev regularity, we prove both existence (as probability densities, essential for the probabilistic interpretation via superposition) and quantitative stability without imposing any growth conditions,
	covering the full range of fractional exponents and extending previous Lagrangian-flow results.
	Second, for drifts whose distributional gradients are given by singular integrals (e.g., Biot--Savart velocity fields from $L^1$ vorticity), we obtain stability estimates that allow the diffusion matrix to be time-space dependent and possibly degenerate, removing the prior restrictions of identity diffusion or bounded coefficients. Our results also recover the pure transport regime as a limiting case.
\end{abstract}



\section{Introduction}\hskip\parindent
We consider the stochastic differential equation (SDE) on $\Rd$:
\begin{equation}\label{SDE}
	\d X_t=b(t,X_t)\,\d t+\sigma(t,X_t)\,\d W_t,
	\quad
	\L_{X_0}=\nu,
\end{equation}
where $\nu\in\P(\Rd)$ and $\L_{X_t}$ denotes the law of $X_t$ for $t\geq0$. Here
$$
b:[0,\infty)\times\Rd\to\Rd,
\qquad
\sigma:[0,\infty)\times\Rd\to\mathbb R^{d\times d}
$$
are measurable, and $W=(W_t)_{t\geq 0}$ is a $d$-dimensional standard Brownian
motion on a complete filtered probability space
$(\Omega,\mathscr F,(\mathscr F_t)_{t\geq 0},\mathbb{P})$. We write $a:=\sigma\sigma^*$
for the diffusion matrix. The Fokker--Planck equation (FPE) associated with \eqref{SDE}
is
\begin{equation}\label{FPE}
	\partial_t\mu_t+\div(b\mu_t)
	=
	\frac12\sum_{i,j=1}^d\partial_{ij}(a_{ij}\mu_t),
	\quad
	\mu_0=\nu .
\end{equation}
It describes how the law of $X_t$ evolves in time, and arises naturally in
probability theory, statistical physics, and fluid dynamics; see, for example,
\cite{Bogachev 2022, Risken 1989}.

When the diffusion term vanishes, that is,
$\sigma=0$, the equation reduces to the continuity equation
\begin{equation}\label{TE}
	\partial_t\mu_t+\div(b\mu_t)=0,
	\quad
	\mu_0=\nu .
\end{equation}
This equation is closely related to the ordinary differential equation
\begin{equation}\label{ODE}
	\d Y_t=b(t,Y_t)\,\d t,
	\quad
	Y_0=x .
\end{equation}
Fix $T>0$. In the classical Cauchy--Lipschitz setting, say
$b\in C^{0,1}_b([0,T]\times \Rd;\R^d)$, the method of characteristics represents the
solution of \eqref{TE} as
$$
\mu_t=\nu\circ Y_t^{-1}, \quad t\in[0,T].
$$
This representation is classical; see, for example, \cite{Crippa 2008}. In many
applications, however, the Lipschitz assumption on $b$ must be weakened. Indeed,
the continuity equation is often part of many equations coming from physics. For
instance, it is part of the inhomogeneous Navier--Stokes equations, where it
expresses conservation of mass. In such models the velocity field is determined
by the system itself and is in general not Lipschitz in space. This motivates the study of \eqref{TE} with rough vector fields, to which we now turn.

The theory begins with DiPerna and Lions \cite{DiPerna 1989}, who proved weak
well-posedness of \eqref{TE} in $L^\infty([0,T];L^q(\Rd))$ for
$b\in L^1([0,T];W^{1,p}(\R^d;\R^d))$ with $1/p+1/q=1$, under suitable bounds on $\div b$
and a mild growth condition.  Ambrosio \cite{Ambrosio 2004} later extended
this theory to vector fields of bounded variation. Both works also established a
qualitative stability result for \eqref{TE}: if a sequence of vector fields $b_n$
converges to $b$ in a suitable sense, then the corresponding solutions $\mu^n$
converge to the solution $\mu$. This stability, however, is qualitative rather
than quantitative. The reason is that the proofs rely on the renormalization
method, which is not suited to deriving quantitative stability estimates.

Crippa and De Lellis \cite{Crippa 2008} obtained the first quantitative stability
estimates for the flow of \eqref{ODE} with a bounded Sobolev vector field $b\in L^\infty([0,T]\times\Rd;\R^d)\cap L^1([0,T];W^{1,p}(\Rd;\R^d))$ for some $p>1$. Their key idea
was to control a logarithmic distance between two flows. More precisely, given two flows $Y^1$ and $Y^2$ generated by vector fields
$b^1$ and $b^2$, they introduced functionals of the form
$$
g_\delta(t)
:=
\left[
\int_{B_r}
\left(
\log\left(
1+\frac{|Y^1(t,x)-Y^2(t,x)|}{\delta}
\right)
\right)^p
\,\d x
\right]^{1/p},
$$
where $p>1$, $r,\delta>0$, and $B_r$ is the open ball in $\R^d$ centered at the origin with radius $r$. With these functionals, they derived stability
estimates of the form
$$
\|Y^1(t,\cdot)-Y^2(t,\cdot)\|_{L^p(B_r)}
\leq
C_t
\left|
\log
\left(
\|b^1-b^2\|_{L^1([0,t];L^p(B_R))}
\right)
\right|^{-1},
$$
where $R=r+T\|b\|_{L^\infty([0,T]\times\Rd)}$, and $C_t$ is a positive constant depending not only on $t$ but also on $r$, $\|b^1\|_{L^\infty([0,t]\times \Rd)}$, $\|b^2\|_{L^\infty([0,t]\times\Rd)}$ and $\|\nabla_x b^1\|_{L^1([0,t];L^p(\Rd))}$. Motivated by this result, Seis \cite{Seis 2017} used ideas from optimal transport to establish quantitative stability estimates for solutions of the continuity equation. Beyond the Sobolev setting, the regularity of $b$ has been relaxed in several directions
\cite{Nguyen 2021, Crippa 2017, Bouchut 2016, Flandoli 2010, Le Bris 2004}; of particular relevance to us, \cite{Bouchut 2013} treated vector fields whose gradients are singular integrals (Assumption \ref{assump stability}). This approach was recently extended by \cite{Crippa 2020} to vector fields with partial fractional Sobolev regularity, though under a stronger restriction on the fractional order than in Assumption \ref{assump1}.

We now turn to the full equation \eqref{FPE}, in which diffusion is present. For
smooth coefficients, the weak well-posedness theory is classical, since linearity
and regularity give existence and uniqueness of weak solutions through standard a
priori estimates \cite{Ladyzenskaja 1968}. The irregular-coefficient case has also been studied
extensively \cite{Bogachev 2022,LeBLions19,Luo 2014,Luo2013,Rockner 2010,Lions 2008,Figalli 2008}.
The coefficients considered in this paper, however,
fall outside the scope of these results. Weak well-posedness for \eqref{FPE} is
therefore established here through the following route. A quantitative stability
estimate gives uniqueness, and an approximation argument  gives existence.

The stability estimate is the main step. We therefore begin
by reviewing the existing quantitative stability results for \eqref{FPE}.
The first named author and Luo \cite{Li 2019} proved such estimates for FPEs with Sobolev
regularity. A key step in their analysis was Trevisan's superposition principle (see Theorem \ref{sp}),
which let them recast the stability problem for the FPE as one for the associated
SDE. They introduced the Kantorovich--Rubinstein distance
$$
\widetilde{\D}_{\delta}(\mu, \nu)
:=
\inf_{\pi \in \mathscr{C}(\mu, \nu)}
\int_{\Rd\times \Rd}
\log\left(1+\frac{|x-y|^2}{\delta^2}\right)
\,\pi(\d x, \d y),
\quad \delta>0,\ \mu, \nu \in \P(\Rd),
$$
where $\mathscr{C}(\mu, \nu)$ denotes the set of all couplings of $\mu$ and $\nu$.
Combining this distance with a synchronous coupling of the two associated SDEs, they obtained
quantitative stability estimates for weak solutions of the FPE. Using the Zvonkin
transformation, they further extended their analysis to singular drifts satisfying
the Ladyzhenskaya--Prodi--Serrin condition. In this setting, the diffusion
coefficient is the identity matrix, that is, $\sigma=I_{d\times d}$, and the
drift coefficient satisfies
$$
b\in L^q([0,T]; L^p(\Rd;\Rd)), \quad \text{with}~p>d,\, q>2~\text{such that}~\frac{d}{p}+\frac{2}{q}<1.
$$
Two other works \cite{Luo 2014, Seis 2022} studied vector fields whose
gradients can be written as a singular integral. Such vector fields were first
studied at the level of the ODE \eqref{ODE} by \cite{Bouchut 2013}, and their analysis was later extended to the FPE by
\cite{Seis 2022} and by   \cite{Luo 2014}. The former requires the diffusion matrix to be the identity, $a=I_{d\times d}$,
while the latter allows a more general $\sigma=\sigma(t,x)$ with Sobolev regularity, at the
cost of an additional boundedness assumption on $b$ and $\sigma$.

The aim of the present paper is to establish the well-posedness and extend the quantitative stability theory for degenerate FPEs in two directions.

First, we treat anisotropic drift vector fields with partial fractional Sobolev regularity of order $s$. For this class, we prove quantitative stability estimates for weak solutions of degenerate FPEs, extending the theory of \cite{Crippa 2020} for Lagrangian flows
of the ODE \eqref{ODE} to the diffusive setting. Compared with that work, our treatment introduces two key improvements: we drop
the growth assumption on the drift $b$, and we cover the full range $s\in(0,1)$, rather than the
more restrictive range $s\in(1/2,1)$ treated there. As a special case, when the diffusion coefficient vanishes
($\sigma=0$), the equation reduces to the continuity equation \eqref{TE}, and our
result gives a quantitative stability estimate for it,  valid for initial data that are
absolutely continuous with respect to the Lebesgue measure.

Second, we study drift vector fields whose distributional derivatives are given by
singular integrals. The typical example is the velocity field reconstructed from an $L^1$ vorticity through the Biot--Savart law,
\begin{equation*}
	b(x)=(k_{\mathrm{BS}}*w)(x), \qquad w\in L^1(\Rd),
\end{equation*}
where the Biot--Savart kernel $k_{\mathrm{BS}}$ is given by
\begin{eqnarray*}
	k_{\mathrm{BS}}(z) =
	\begin{cases}
		\dfrac{z^{\perp}}{2\pi |z|^2}, & d = 2,\ z \ne 0,\\[3ex]
		\dfrac{z}{s_{d-1}|z|^{d}}, & d \ge 3,\ z \ne 0,
	\end{cases}
\end{eqnarray*}
with $z^{\perp}:=(-z_2,z_1)$ for $z=(z_1,z_2)\in \R^2$, and $s_{d-1}$ the area of the
$(d-1)$-dimensional unit sphere. This kernel is of independent interest and plays a central role in fluid dynamics, as it recovers the velocity field from the vorticity; through this relation, the incompressible Navier--Stokes equation can be written in vorticity form. Drift fields of this type therefore arise naturally in applications. For this class, we again establish quantitative stability estimates for weak solutions of degenerate FPEs. Our result improves on \cite{Seis 2022, Luo 2014} in two respects. First, we allow a
time- and space-dependent diffusion matrix $\sigma=\sigma(t,x)$ that may be degenerate, rather than only the identity matrix required in previous works. Second, we remove the boundedness assumption on the coefficients that was imposed in \cite{Luo 2014}. In addition, as before, our estimate reduces to that for the continuity  equation \eqref{TE}  when $\sigma=0$.

In both settings, the stability estimate directly implies uniqueness within the corresponding class of weak solutions. Existence follows from the DiPerna--Lions theory, which was first introduced by DiPerna and Lions in their study of the Boltzmann equation \cite{DiPernaLions1989 boltzmann, DiPernaLions1988 boltzmann}. This theory yields a weak solution, but it does not ensure that the solution  is a probability measure. We establish this additional property, under the mild assumption of a finite first moment for the initial datum. The point is irrelevant in the classical PDE theory but essential here, since the superposition principle applies precisely to probability-measure-valued solutions. Through it, a stability estimate for \eqref{FPE} becomes one for the associated SDE \eqref{SDE}. At the level of the SDE, the estimate can be obtained by probabilistic tools such as coupling, which do not require the diffusion matrix to be the identity. This is in contrast to the analytic argument of \cite{Seis 2022}, which relies on that assumption, and it allows us to cover a broader, possibly degenerate class of equations. Our results also complement \cite{Luo 2014}, which established uniqueness of weak solutions but not their existence. Together, these results give the weak well-posedness of the degenerate FPEs considered in this paper.

\textbf{Structure of the paper}\hskip\parindent
This paper is organized in the following way. In Section \ref{pre}, we collect the necessary preliminaries by recalling basic definitions and presenting some useful results that will be used in the following proofs. Section \ref{par} is devoted to the first class of drift coefficients, which have partial fractional Sobolev regularity. The main result (Theorem \ref{thm partial}) establishes both existence of a weak solution to the degenerate FPE and a quantitative stability estimate for two such solutions. The proof of the existence part relies on a compactness argument via the DiPerna--Lions theory, while the stability estimate is obtained by coupling the associated SDEs and using anisotropic Kantorovich--Rubinstein distances. At the end of the section, in Subsection \ref{rate}, we discuss the rates of convergence that follow from the stability estimate. In Section \ref{singular}, we treat the second class of drifts, whose distributional gradients are given by singular integrals. The main theorem (Theorem \ref{stability estimate 2}) provides a quantitative stability estimate for weak solutions under these rough coefficients. Its proof follows a similar SDE-coupling strategy, but requires a refined maximal function estimate adapted to singular kernels. Finally, in Subsection \ref{app}, we apply the stability result to prove well-posedness of the FPE under an additional structural assumption on the drift, covering the physically relevant case of velocity fields generated from an $L^1$ vorticity via the Biot--Savart law.

\textbf{Notation}\hskip\parindent
Throughout the paper, $d\geq2$ is a fixed positive integer, and $\mathbb{R}^d$ denotes the $d$-dimensional Euclidean space with inner product $\langle\cdot,\cdot\rangle$ and norm $|\cdot|$, and $\mathcal{L}^d$ is the Lebesgue measure on $\mathbb{R}^d$. For $x\in\R^d$ and $R>0$, $B_R(x):=\{y\in\mathbb{R}^d: |y-x|<R\}$, and write $B_R:=B_R(0)$. We use the following conventions and function spaces.

\medskip

\noindent\textbf{Basic notation}
\begin{itemize}[itemsep=0pt]
	\item[$\cdot$] For a set $A$, $\mathbbm{1}_A$ denotes its  indicator function.
	
	\item[$\cdot$] For a Polish space $M$, $\mathcal{P}(M)$ is the set of Borel probability measures on $M$, equipped with the weak topology. 
	\item[$\cdot$] For a matrix $Q$, denote by $Q^*$ its transpose, by $\|Q\|$ its Hilbert--Schmidt norm,  by $Q_{ij}$ its entries,  and (if $Q$ is square) by $\operatorname{tr}Q$ its trace.
	\item[$\cdot$] $f\lesssim g$ means $f\le C g$ for some constant $C>0$ depending only on fixed parameters.
	\item[$\cdot$] $f*g$ is the convolution of $f$ and $g$. In particular, for $n=1,2,\cdots$, if $f$ is $\mathbb{R}^n$-valued and $g$ is $\mathbb{R}$-valued, then $f*g=g*f=(f_1*g,\dots,f_n*g)$.
	\item[$\cdot$] For a function $f$, $f^-:=-\min\{f,0\}$.
	\item[$\cdot$] For a random variable $Z$, $\mathscr{L}_{Z}$ denotes the law of $Z$.
	\item[$\cdot$] $\nabla,\mathrm{div},\nabla^2$ and $\partial_j$ denote the spatial gradient, divergence, Hessian, and $j$-th partial derivative, respectively.
\end{itemize}

\noindent\textbf{Scalar function spaces on \(\mathbb{R}^d\)}
\begin{itemize}[itemsep=0pt]
	\item[$\cdot$] $L^p(\mathbb{R}^d)$ (\(1\le p\le\infty\)) denotes the standard Lebesgue space, with norm $\|\cdot\|_{L^p(\R^d)}$. $L^p_{\mathrm{loc}}(\R^d)$ denotes the local version.
	\item[$\cdot$] $W^{k,p}(\mathbb{R}^d)$ ($k=1,2,\cdots$, $1\le p\le\infty$) denotes the Sobolev space of functions whose weak derivatives up to order $k$ lie in $L^p(\R^d)$, with norm $\|\cdot\|_{W^{k,p}(\R^d)}$. $W^{k,p}_{\mathrm{loc}}(\R^d)$ denotes the local version.
	\item[$\cdot$] For $0<s<1$, $1\le p<\infty$, the fractional Sobolev space $W^{s,p}(\mathbb{R}^d)$ is defined by
	\begin{align*}
		&W^{s,p}(\R^d):=\big\{f\in L^p(\R^d):\ [f]_{W^{s,p}(\R^d)}<\infty\big\},\\
		&[f]_{W^{s,p}(\Rd)}^p:=\int_{\mathbb{R}^d}\int_{\mathbb{R}^d}\frac{|f(x)-f(y)|^p}{|x-y|^{d+sp}}\,\d x \, \d y,
	\end{align*}
	with norm $\|f\|_{W^{s,p}(\R^d)}:=\|f\|_{L^p(\R^d)}+[f]_{W^{s,p}(\Rd)}$.
	\item[$\cdot$] $\mathscr{S}(\R^d)$ denotes the Schwartz space of rapidly decreasing smooth functions, and its dual $\mathscr{S}'(\R^d)$ denotes the space of tempered distributions.
	\item[$\cdot$] $\mathscr{D}(\Rd)=C_c^\infty(\Rd)$ denotes the space of compactly supported smooth functions, and its dual $\mathscr{D}'(\Rd)$ denotes the space of distributions on $\Rd$.
	\item[$\cdot$] $C_b(\mathbb{R}^d)$ denotes the class of bounded continuous functions. $C_b^k(\mathbb{R}^d)$ ($k=1,2,\cdots$) denotes the class of continuous functions bounded with bounded derivatives up to order $k$. $C_b^\infty(\R^d)$ denotes the class of continuous functions bounded with all derivatives (not a Banach space). $C^\infty(\mathbb{R}^d)$ denotes the class of smooth functions.
\end{itemize}

\noindent\textbf{Time-space function spaces}\hskip\parindent
Let $J\subset\mathbb{R}$ be an interval.
\begin{itemize}[itemsep=0pt]
	\item[$\cdot$] $L^p(J\times\mathbb{R}^d)$ denotes the Lebesgue space with respect to the product Lebesgue measure $\mathcal{L}^1\otimes\mathcal{L}^d$.  $L^p_{\mathrm{loc}}(J\times\mathbb{R}^d)$ is the local version.
	\item[$\cdot$]  $C^{0,1}_b(J\times \Rd)$ denotes bounded continuous functions together with their first space derivatives.  $C_b^{1,2}(J\times\mathbb{R}^d)$ denotes bounded continuous functions together with their first time derivative and second space derivatives. $C_c^{1,2}(J\times\mathbb{R}^d)$ denotes the subspace with compact support in the space variable (and vanishing near the time endpoint if $J=[0,T)$).
	\item[$\cdot$] (\textit{Bochner spaces}) For a Banach space $(E,\|\cdot\|_E)$, the space $L^q(J;E)$ ($1\le q\le\infty$) consists of strongly measurable $f:J\to E$ with
	\begin{align*}
		\|f\|_{L^q(J;E)}:=
		\begin{cases}
			\left(\int_J \|f(t)\|_E^q\,dt\right)^{1/q}, & 1\le q<\infty,\\[4pt]
			\operatorname{ess\,sup}_{t\in J}\|f(t)\|_E, & q=\infty.
		\end{cases}
	\end{align*}
	We use the shorthand (mainly in proofs)
	\begin{align*}
		L^q_t L^p_x:=L^q(J;L^p(\mathbb{R}^d)),\quad
		L^q_t W^{k,p}_x:=L^q(J;W^{k,p}(\mathbb{R}^d)),\quad
		L^q_t W^{s,p}_x:=L^q(J;W^{s,p}(\mathbb{R}^d)),
	\end{align*}
	and the corresponding norms
	\begin{align*} \|\cdot\|_{L^q_tL^p_x},\quad
		\|\cdot\|_{L^q_tW^{k,p}_x},\quad
		\|\cdot\|_{L^q_tW^{s,p}_x}.
	\end{align*}
	Analogous notation applies when $E$ is replaced by any other space.
	\item[$\cdot$] (\textit{Convention for \(C_b^\infty\) in Bochner spaces}) The notation $L^q(J;C_b^\infty(\mathbb{R}^d))$ is an abuse, meaning that for a.e. $t\in J$, $f(t,\cdot)\in C_b^\infty(\R^d)$, and the $L^q$-norm of its derivatives up to the required finite order (which is always fixed by the context) is integrable on $J$.
\end{itemize}

\medskip

\noindent\textbf{Weak Lebesgue / Lorentz--Marcinkiewicz spaces}\hskip\parindent
For $1\le r<\infty$ and a measure space $(\mathcal{X},\mu)$, $L^{r,\infty}(\mathcal{X},\mu)$ denotes the space of all measurable functions $f:\mathcal{X}\to \R$ such that
$$
\|f\|_{L^{r,\infty}(\mathcal{X},\mu)}:=\sup_{\lambda>0}\lambda\,\mu(\{|f|>\lambda\})^{1/r}<\infty.
$$
For $r=\infty$, write $L^{\infty,\infty}:=L^\infty$. On $\mathbb{R}^d$, we write $L^{r,\infty}(\mathbb{R}^d)$, and time-dependent versions are given by the Bochner construction above, e.g., $L^q(J;L^{r,\infty}(\mathbb{R}^d))$ with
\begin{align*}
	\|f\|_{L^q(J;L^{r,\infty}(\mathcal{X},\mu))}
	:=\begin{cases}
		\left(\int_J\|f(t,\cdot)\|_{L^{r,\infty}(\mathcal{X},\mu)}^q\,\d t\right)^{1/q}, & 1\le q<\infty,\\[4pt]
		\operatorname{ess\,sup}_{t\in J}\|f(t,\cdot)\|_{L^{r,\infty}(\mathcal{X},\mu)}, & q=\infty.
	\end{cases}
\end{align*}

\medskip

\noindent\textbf{Vector- and matrix-valued analogues}\hskip\parindent
All the spaces above have natural componentwise extensions to functions taking values in $\mathbb{R}^m$ or $\mathbb{R}^{m\times m}$ ($m=1,2,\cdots$). Norms are taken as the Euclidean norm (for vectors) or the Hilbert--Schmidt norm (for matrices). We often omit the target space when it is clear from the context (e.g., $b$ is vector-valued, $\sigma$ is matrix-valued).

\section{Preliminaries}\label{pre}\hskip\parindent
We begin by recalling the definitions of weak solution for the FPE and for the SDE, respectively.  For $f\in C^{1,2}((0,T)\times\R^d)$, define the operator $L$ by
\begin{align*}
	( Lf)(t,x)=
	\frac{1}{2}\operatorname{tr}\left((\sigma\sigma^*)(t,x)[\nabla^2 f(t,\cdot)](x)\right)
	+
	\< b(t,x),[\nabla f(t,\cdot)](x)\>,\quad t\in(0,T),\,x\in\R^d.
\end{align*}
\begin{definition}[Weak solution of the FPE]\label{def weak solution FPE}
	Let $T>0$ and $\nu\in\P(\Rd)$. A Borel curve
	$(\mu_t)_{t\in[0,T]}\subset\P(\Rd)$ is called a weak solution to \eqref{FPE}
	if the following conditions hold:
	\begin{itemize}
		\item[(i)](Integrability condition)
		\begin{equation}\label{integral condition}
			\int_0^T\int_{\Rd}
			\left(\|a(t,x)\|+|b(t,x)|\right)
			\mu_t(\d x)\,\d t
			<\infty;
		\end{equation}

		\item[(ii)](Weak formulation) For every $\varphi\in C_c^{1,2}([0,T)\times\Rd)$,
		\begin{equation}\label{weak solution integral}
			\int_{\Rd} \varphi(0,x)\,\nu(\d x)
			+
			\int_0^T\int_{\Rd}
			\left(\partial_t \varphi(t,x)+L\varphi(t,x)\right)
			\mu_t(\d x)\,\d t
			=0.
		\end{equation}
	\end{itemize}
\end{definition}

\begin{remark}\label{narrowly continuous}
	Recall that a curve $(\mu_t)_{t\in[0,T]}\subset\P(\Rd)$ is called weakly
	continuous if for every $g\in C_b(\Rd)$, the map
	\[
	t\mapsto \int_{\Rd} g(x)\,\mu_t(\d x)
	\]
	is continuous. By a standard density argument, see, e.g.,
	\cite[Lemma~8.1.2]{Ambrosio 2008} and \cite[Remark~2.3]{Trevisan 2016}, every
	weak solution admits a unique weakly
	continuous representative
	$(\tilde\mu_t)_{t\in[0,T]}\subset\P(\Rd)$, with $\tilde\mu_t=\mu_t$ for a.e.
	$t\in[0,T]$ and $\tilde\mu_0=\nu$. Hence, it is meaningful to say that the solution starts from $\nu$.
\end{remark}

If $(\mu_t)_{t\in[0,T]}$ is absolutely continuous with respect to the Lebesgue
measure, we write
\[
\mu_t(\d x)=\rho(t,x)\,\d x .
\]
In this case we also call $\rho$ a weak solution to the FPE \eqref{FPE}.

\begin{definition}[Weak solution of the SDE]\label{def weak solution SDE}
	Let $T>0$ and $\nu\in\P(\Rd)$. A pair $(X,W)$ is called a weak solution
	to the SDE \eqref{SDE} with initial distribution $\nu$ if there exists a complete
	filtered probability space
	$
	(\Omega,\mathscr F,(\mathscr F_t)_{t\in[0,T]},\mathbb P)
	$
	such that $W=(W_t)_{t\in[0,T]}$ is a $d$-dimensional standard $(\mathscr F_t)_{t\in[0,T]}$-Brownian motion,
	$X=(X_t)_{t\in[0,T]}$ is an $(\mathscr F_t)_{t\in[0,T]}$-adapted continuous
	process, $\L_{X_0}=\nu$, the integrability condition
	\[
	\E\int_0^T
	\left(|b(s,X_s)|+\|\sigma(s,X_s)\|^2\right)\,\d s<\infty
	\]
	holds, and, for every $t\in[0,T]$, $\mathbb P$-a.s.,
	\[
	X_t
	=
	X_0
	+
	\int_0^t b(s,X_s)\,\d s
	+
	\int_0^t \sigma(s,X_s)\,\d W_s .
	\]
	For brevity, we shall often write simply $(X,W)$. When the underlying probability space needs to be specified explicitly, we use the full notation $(\Omega,\mathscr F,(\mathscr F_t)_{t\in[0,T]},\mathbb P;X,W)$.
\end{definition}

The link between these two notions of solution is provided by the martingale
solutions of the SDE, which we now define. Let $C([0,T];\Rd)$ be the space of
continuous paths from $[0,T]$ to $\Rd$. For $0\leq t\leq T$, let
$$
\mathcal B_t(C([0,T];\Rd))
:=
\sigma(e_s:0\leq s\leq t)
$$
be the canonical $\sigma$-algebra, where
$e_s: C([0,T];\Rd)\rightarrow\R^d$ is the evaluation map at time $s\in[0,T]$.

\begin{definition}[Martingale solution of the SDE]\label{def martingale solution SDE}
	Given $\nu\in\P(\Rd)$, a probability measure
	$Q\in\P(C([0,T];\Rd))$ is called a martingale solution to the SDE
	\eqref{SDE} with initial distribution $\nu$ if
	$$Q\circ e_0^{-1}=\nu,$$
	$$
	\int_0^T \int_{\Rd}
	\left(|b(t,x)|+\|a(t,x)\|\right)
	(Q\circ e_t^{-1})(\d x)\,\d t
	<\infty,
	$$
	and, for every $\varphi\in C_b^{1,2}([0,T]\times\Rd)$, the process
	$(M_t^\varphi)_{t\in[0,T]}$ defined by
	$$
	M_t^\varphi(\omega)
	:=
	\varphi(t,\omega_t)-\varphi(0,\omega_0)
	-
	\int_0^t(\partial_s \varphi+ L\varphi)(s,\omega_s)\,\d s,\quad \omega\in C([0,T];\Rd),
	$$
	is a $\mathcal B_t(C([0,T];\Rd))$-martingale under $Q$.
\end{definition}

It is well known that martingale solutions and weak solutions of the SDE are
equivalent in the following sense. If $Q\in\P(C([0,T];\Rd))$ is a martingale
solution to \eqref{SDE}, then there exists a weak solution
$(\Omega,\mathscr F,(\mathscr F_t)_{t\in[0,T]},\mathbb P; X,W)$ to \eqref{SDE}
with $\mathbb P\circ X^{-1}=Q$. Conversely, if $(X,W)$ is a weak solution to
\eqref{SDE}, then its path law $\mathbb P\circ X^{-1}$ is a martingale solution.
For a comprehensive treatment, we refer to the monographs \cite{StroockVaradhan 2006,Ikeda 1989}.

On the one hand, for any $\varphi\in C_c^{1,2}([0,T)\times\Rd)$, applying It\^o's
formula to $\varphi(t,X_t)$ gives \eqref{weak solution integral}. Hence the time
marginals of any weak solution to the SDE form a weak solution to the
corresponding FPE \eqref{FPE}. On the other hand, the converse direction follows
from the following superposition principle.

\begin{theorem}[{Superposition principle, \cite[Theorem~2.5]{Trevisan 2016}}]\label{sp}
	Let $(\mu_t)_{t\in[0,T]}\subset\P(\Rd)$ be a weakly continuous weak solution to the FPE \eqref{FPE}. Then there exists a martingale solution
	$Q\in\P(C([0,T];\Rd))$ to the SDE \eqref{SDE} such that
	$$Q\circ e_t^{-1}=\mu_t,\qquad  t\in[0,T].$$
	Equivalently, there exists a weak solution $(X,W)$ to the SDE satisfying
	$$\L_{X_t}=\mu_t,\qquad t\in[0,T].$$
\end{theorem}

Our stability estimates are formulated in terms of expectations of the weak solutions to the associated SDEs. These weak solutions may, however, be defined on different probability spaces, and hence cannot be compared directly. To address this issue, we use the following lemma (cf. \cite[Proposition 3.6]{Li 2019}). It constructs two processes on a common filtered probability space, driven by the same Brownian motion, each having the same law as one of the original weak solutions. This coupling serves as the starting point of our stability analysis.
\begin{lemma}\label{lem-common}
	Let
	$(\Omega,\mathscr F,(\mathscr F_t)_{t\geq 0},\mathbb P;X,W)$
	and
	$(\hat{\Omega},\hat{\mathscr F},(\hat{\mathscr F}_t)_{t\geq 0},\hat{\mathbb P};\hat X,\hat W)$
	be weak solutions to the SDE \eqref{SDE} corresponding to the coefficients
	$(b,\sigma)$ and $(\hat b,\hat\sigma)$, with initial distributions
	$\nu$ and $\hat\nu$, respectively. Then, for any
	$\pi_0\in\mathscr C(\nu,\hat\nu)$, there exist a filtered probability space
	$
	(\tilde{\Omega},\tilde{\mathscr F},
	(\tilde{\mathscr F}_t)_{t\geq 0},\tilde{\mathbb P}),
	$
	a standard $d$-dimensional Brownian motion $\tilde W=(\tilde W_t)_{t\geq0}$ adapted to
	$(\tilde{\mathscr F}_t)_{t\geq0}$, and two
	$(\tilde{\mathscr F}_t)_{t\geq0}$-adapted continuous processes
	$Y$ and $\hat Y$ such that
	$\L_{(Y_0,\hat Y_0)}=\pi_0$ and the following properties hold:
	\begin{itemize}
		\item[(i)] The pairs $(X,W)$ and $(Y,\tilde W)$ have the same law, and so do
		the pairs $(\hat X,\hat W)$ and $(\hat Y,\tilde W)$.
		
		\item[(ii)] The systems
		$
		(\tilde{\Omega},\tilde{\mathscr F},
		(\tilde{\mathscr F}_t)_{t\geq 0},\tilde{\mathbb P};Y,\tilde W)
		$
		and
		$
		(\tilde{\Omega},\tilde{\mathscr F},
		(\tilde{\mathscr F}_t)_{t\geq 0},\tilde{\mathbb P};\hat Y,\tilde W)
		$
		are weak solutions to the SDE \eqref{SDE} with coefficients
		$(b,\sigma)$ and $(\hat b,\hat\sigma)$, and with initial distributions
		$\nu$ and $\hat\nu$, respectively.
	\end{itemize}
\end{lemma}

In view of this lemma, we may assume without loss of generality that the weak solutions under consideration are defined on a common filtered probability space and are driven by the same Brownian motion.

Finally, let $(\mathcal{X},\mu)$ be a measure space. For $r\in[1,\infty]$, we denote by $L^r(\mathcal{X},\mu)$ the usual Lebesgue space equipped with the norm $\|\cdot\|_{L^r(\mathcal{X},\mu)}$, and by $L^{r,\infty}(\mathcal{X},\mu)$ the weak Lebesgue / Lorentz--Marcinkiewicz space endowed with the quasi-norm $\|\cdot\|_{L^{r,\infty}(\mathcal{X},\mu)}$; precise definitions are given in \textbf{Notation}. Recall the continuous embedding $L^r(\mathcal{X},\mu)\subset L^{r,\infty}(\mathcal{X},\mu)$, which satisfies
$$\|f\|_{L^{r,\infty}(\mathcal{X},\mu)}\leq\|f\|_{L^r(\mathcal{X},\mu)},\quad f\in L^r(\mathcal{X},\mu).$$
We shall also make use of the following interpolation inequalities; for their proofs, we refer to \cite{Crippa 2017,Bouchut 2013}.
\begin{lemma}[Interpolation]\label{Interpolation}
	Let $E\subset\mathcal{X}$ be a measurable set with $\mu(E)<\infty$,  and let $f:\mathcal{X}\rightarrow\R$ be a measurable function. Then the following estimates hold.
	\begin{itemize}
		\item[(1)] For $1\leq r<\ell\leq\infty$,
		\begin{equation}\label{interp 1}
			\|f\|_{L^r(E,\mu)}^r
			\leq
			\frac{\ell}{\ell-r}
			\mu(E)^{1-\frac{r}{\ell}}
			\|f\|_{L^{\ell,\infty}(E,\mu)}^r.
		\end{equation}

		\item[(2)] For $1<\ell\leq\infty$,
		\begin{equation}\label{interp 2}
			\|f\|_{L^1(E,\mu)}
			\leq
			\frac{\ell}{\ell-1}
			\|f\|_{L^{1,\infty}(E,\mu)}
			\left[
			1+
			\log\left(
			\frac{
				\mu(E)^{1-\frac{1}{\ell}}
				\|f\|_{L^{\ell,\infty}(E,\mu)}
			}{
				\|f\|_{L^{1,\infty}(E,\mu)}
			}
			\right)
			\right].
		\end{equation}
	\end{itemize}
	Here, for $\ell=\infty$, the factors $\frac{\ell}{\ell-r}$ in \eqref{interp 1} and $\frac{\ell}{\ell-1}$ in \eqref{interp 2} are interpreted as their limiting values, both equal to $1$, and if $\|f\|_{L^{1,\infty}(E,\mu)}=0$, the right-hand side of \eqref{interp 2} is understood as $0$ by continuity.
\end{lemma}

Further details on weak Lebesgue spaces can be found in the monograph \cite{Grafakos2014a}.

\section{Partially regular drift coefficients}\label{par}
\subsection{Main results}\hskip\parindent
In this subsection, we present our main results for FPEs whose drift coefficients have partial fractional Sobolev regularity. To this end, we split the state space as
$$\Rd=\R^{d_1}\times\R^{d_2}, \quad d=d_1+d_2,$$
where $d_1,d_2$ are positive integers. Consider drift coefficients $b:[0,T]\times\R^{d_1}\times\R^{d_2}\rightarrow \R^{d_1}\times\R^{d_2}$ of the form
$$
b(t,x_1,x_2)
=
\bigl(b_1(t,x_1), b_2(t,x_1,x_2)\bigr),
$$
and assume that the diffusion coefficient $\sigma: [0,T]\times\R^{d_1}\times\R^{d_2}\rightarrow \R^{d\times d}$ has the block-diagonal form
$$
\sigma(t,x_1,x_2)
=
\begin{pmatrix}
	\underbrace{\sigma_1(t,x_1)}_{d_1\times d_1} & 0 \\[5ex]
	0 & \underbrace{\sigma_2(t,x_1,x_2)}_{d_2\times d_2}
\end{pmatrix},
$$
where $t\in[0,T]$ and $(x_1,x_2)\in\R^{d_1}\times\R^{d_2}$.
Note that $b_1$ and $\sigma_1$ do not depend on  the second variable $x_2$.

Under this structural assumption, we impose the following conditions on the
coefficients. For convenience, we define the modified drift
$$b^\sigma:=b-\frac{1}{2}\div(\sigma\sigma^*).$$

\begin{assumption}\label{assump1}
	Let $p,r\in (1,\infty]$, $s\in(0,1)$, and $\nu\in\mathcal{P}(\R^d)$.
	\begin{itemize}
		\item[]$(\mathcal{H}^{\sigma})$
		$\sigma \in L^2([0,T]; L^\infty(\R^d;\R^{d\times d}))$ with
		\begin{align*}
			&\bigg( \sum_{i,j,k=1}^{d_1} |\partial_k (\sigma_1)_{ij}|^2 \bigg)^{1/2}\in L^2([0,T]; L^{2r}(\R^{d_1}) ) ,\\
			&\bigg( \sum_{k=1}^d\sum_{i,j=1}^{d_2}|\partial_k (\sigma_2)_{ij}|^2 \bigg)^{1/2}\in L^2([0,T];  L^{2r}(\R^d)),
		\end{align*}
		where the	derivative is understood in the weak sense.
		
		\item[]$(\mathcal{H}^{b})$
		$b \in L^p_{\mathrm{loc}}([0,T]\times\mathbb{R}^d;\R^d)$ with
		\begin{align*}
			b_1 &\in  L^1([0,T]; W^{1,p}(\R^{d_1}; \R^{d_1})) \cap L^1([0,T]; L^\infty(\R^{d_1};\R^{d_1})),\\
			b_2 & \in L^1([0,T]\times \R^{d_2}; W^{s,1}(\R^{d_1}; \R^{d_2}))\cap L^1([0,T]\times \R^{d_1}; W^{1,p}(\R^{d_2};\R^{d_2}) ). 	
		\end{align*}
		
		\item[]$(\mathcal{H}^{b^\sigma})$
		$(\operatorname{div} b^\sigma)^- \in L^1([0,T]; L^\infty(\R^d))$.
		
		\item[]$(\mathcal{H}^{\nu})$
		The initial distribution $\nu$ admits a probability density
		$\rho_0 \in L^1(\R^d) \cap L^\infty(\R^d)$ with
		$$\int_{\mathbb{R}^d}|x|\rho_0(x)\,\d x < \infty.$$
	\end{itemize}
\end{assumption}

\begin{remark}\label{rem-assump1}
	(i) The condition $(\div b)^-\in  L^1([0,T]; L^\infty(\R^d))$ is standard in the study of
	transport equations; see, for instance, \cite{Ambrosio 2004,DiPerna 1989}. It has
	a natural interpretation in terms of the flow generated by $b$: the divergence
	quantifies
	local compression and expansion. The most familiar case is $\div b=0$,
	which corresponds to an incompressible flow; the condition on $(\div b)^-$ is
	weaker, allowing compression but keeping it controlled. In the FPE setting, as
	shown by Le Bris and Lions \cite{Lions 2008}, the same idea leads naturally to
	condition $(\mathcal H^{b^\sigma})$. This assumption provides the
	fundamental a priori estimate for densities, which is essential for constructing
	weak solutions when the drift has low regularity.

	(ii) The conditions $\int_{\Rd}|x|\rho_0(x)\,\d x<\infty$ and
	\eqref{integral condition} together imply that the solution to the FPE has a
	finite first moment, namely
	$$
	\operatorname*{ess\,sup}_{t\in[0,T]}\int_{\Rd}|x|\mu_t(\d x)<\infty.
	$$
	Indeed, by the superposition principle, let $(X,W)$ be a weak solution to the
	corresponding SDE with $\L_{X_t}=\mu_t$ for a.e. $t\in[0,T]$. Condition
	\eqref{integral condition} guarantees that $\int_0^t\sigma(s,X_s)\,\d W_s$ is a
	square-integrable martingale; see, for instance, \cite{Ikeda 1989}. By Jensen's
	inequality and It\^o's isometry, we obtain, for a.e. $t\in[0,T]$,
	\begin{equation*}\label{SDE moment}
		\begin{split}
			\int_{\Rd}|x|\mu_t(\d x) &= \E |X_t|\\
			&\leq \E|X_0|
			+\E\left|\int_0^t\sigma(s,X_s)\,\d W_s\right|
			+\E\left|\int_0^t b(s,X_s)\,\d s\right|\\
			&\leq \int_{\Rd}|x|\rho_0(x)\,\d x
			+\left(\int_0^T\int_{\Rd}\|\sigma(s,x)\|^2\mu_s(\d x)\,\d s\right)^{\frac12}\\
			&\quad +\int_0^T\int_{\Rd}|b(s,x)|\, \mu_s(\d x)\,\d s\\
			&<\infty .
		\end{split}
	\end{equation*}
	Taking the essential supremum in time gives the desired estimate.

	(iii) Under assumption $(\mathcal{H}^{\sigma})$, one readily verifies that for a.e. $t\in[0,T]$,
	$$\sigma_1(t,\cdot)\in W^{1,1}_{\rm loc}(\R^{d_1}; \R^{d_1\times d_1}),\qquad\sigma_2(t,\cdot)\in W^{1,1}_{\rm loc}(\R^{d}; \R^{d_2\times d_2}).$$  This regularity makes Lemma \ref{Morrey-type estimate} applicable.
\end{remark}

Corresponding to the structure of the coefficients $b,\sigma$, for $\delta_1, \delta_2>0$, we introduce the
$d\times d$ diagonal matrix
\begin{align}\label{diag-B}
	B = \mathrm{diag}(\underbrace{\delta_1, \dots, \delta_1}_{d_1\text{ times}}, \underbrace{\delta_2, \dots, \delta_2}_{d_2\text{ times}}),
\end{align}
which defines an anisotropic Kantorovich--Rubinstein distance
$$\D_{\delta_1, \delta_2}(\mu, \nu):=\inf_{\pi \in \mathscr{C}(\mu, \nu)}\int_{\Rd\times \Rd} \log\left(1+|B^{-1}(x-y)|\right) \,\pi(\d x, \d y), \quad \mu, \nu \in \P(\Rd),$$
where $B^{-1}$ denotes the inverse matrix of $B$.
When the measures admit densities $\rho_1,\rho_2$, we also write
$\D_{\delta_1, \delta_2}(\rho_1,\rho_2)$ for the distance between the
corresponding distributions.

We now present the main theorem of this section.
It has two components: first, an existence result for weak solutions that are probability densities, a property crucial for invoking the superposition principle and thereby linking the FPE to the associated SDE; second, a quantitative stability estimate for any two such solutions. Owing to the local integrability of the coefficients $b_1$ and $\sigma_1$, the latter is formulated in terms of expectations of the corresponding SDE solutions, with trajectories stopped upon exiting a ball in the state space.
\begin{theorem}\label{thm partial}
	The following statements hold.
	\begin{itemize}
		\item[(1)] \textbf{(Existence)}
		Under Assumption \ref{assump1}, the FPE \eqref{FPE} admits a weak solution that is absolutely continuous with respect to the Lebesgue measure, with a probability density $\rho\in L^\infty([0,T]; L^1(\R^d)\cap L^\infty(\R^d))$.
		
		\item[(2)] \textbf{(Quantitative stability)}	
		Let $(b,\sigma,\rho_0)$ and $(\hat b,\hat\sigma,\hat\rho_0)$ be two triples satisfying Assumption \ref{assump1}. Denote by $\rho,\hat\rho\in L^\infty([0,T]; L^1(\R^d)\cap L^\infty(\R^d))$ the corresponding weak solutions to the FPE \eqref{FPE}, and by $(X,W)$ and $(\hat X,W)$ the associated weak solutions to the SDE \eqref{SDE}. For $R>0$, define the stopping time
		$$
		\tau_R:=\inf\{t\in[0,T]:\ |X_t|\vee|\hat X_t|\ge R\}.
		$$
		Then, for every $\delta_2\in(0,1)$, one can choose
		$\delta_1\in(0,\delta_2)$ such that, for all $R>0$, the following stability
		estimate holds:
		\begin{equation*}
			\begin{split}
				&\sup_{t\in[0,T]}
				\E\left\{\log\left(1+
				\left|B^{-1}(X_{t\wedge \tau_R}-\hat X_{t\wedge \tau_R})\right|
				\right)\right\}  \\
				&\quad \lesssim
				\D_{\delta_1,\delta_2}(\rho_0,\hat\rho_0)
				+
				\frac{1}{\delta_1}
				\|b_1-\hat b_1\|_{L^1([0,T];L^p(B_R))}  \\
				&\qquad
				+
				\frac{1}{\delta_1}
				\|b_2-\hat b_2\|_{L^1([0,T];L^1(\Rd))}
				+
				\frac{1}{\delta_1^2}
				\|\sigma-\hat\sigma\|^2_{L^2([0,T]; L^{2r}(B_R))} +1.
			\end{split}
		\end{equation*}
		Here the implicit constant in \(\lesssim\) is independent of $\delta_1,\delta_2$, and depends only on \(d,T,R\), the $L^\infty([0,T]; L^1(\R^d)\cap L^\infty(\R^d))$-norms of $\rho,\hat\rho$, and
		\[
		\begin{aligned}
			&\|b\|_{L^p([0,T]\times B_{R+1})},\quad
			\|b_2\|_{L^1([0,T]\times\mathbb{R}^{d_2};W^{s,1}(\R^{d_1}))},\quad
			\|\nabla_{x_1}b_1\|_{L^1([0,T];L^p(\R^{d_1}))},
			\\
			&\|\nabla_{x_2}b_2\|_{L^1([0,T]\times\mathbb{R}^{d_1};L^p(\R^{d_2}))},\quad
			\|\nabla_{x_1}\sigma_1\|_{L^2([0,T];L^{2r}(\R^{d_1}))},\quad
			\|\nabla_x\sigma_2\|_{L^2([0,T];L^{2r}(\R^{d}))}.
		\end{aligned}
		\]
	\end{itemize}
\end{theorem}

\begin{remark}
	\begin{itemize}
		\item[(1)] The above stability estimate implies uniqueness of weak solutions to the FPE \eqref{FPE} in the class $L^\infty([0,T]; L^1(\R^d)\cap L^\infty(\R^d))$. Indeed,
		if $b=\hat b$, $\sigma=\hat\sigma$, and $\rho_0=\hat\rho_0$, then the right-hand side of the estimate is bounded independently of
		$\delta_1$ and $\delta_2$. Letting $\delta_1,\delta_2\downarrow0$, the
		left-hand side would diverge unless
		$X_{t\wedge\tau_R}=\hat X_{t\wedge\tau_R}$ a.s. for every $t\in[0,T]$.
		Letting $R\to\infty$ gives $X_t=\hat X_t$ a.s. for every
		$t\in[0,T]$, and hence $\rho=\hat\rho$ in
		$L^\infty([0,T]; L^1(\R^d)\cap L^\infty(\R^d))$.
		
		\item[(2)] The assumption $b_1 \in L^1([0,T]; L^\infty(\R^{d_1};\R^{d_1}))$ in $(\mathcal{H}^{b})$ is used only to
		guarantee the integrability condition \eqref{integral condition}. It may be
		replaced by any other assumption that ensures the same condition.
	\end{itemize}
\end{remark}

\subsection{Proof of Theorem \ref{thm partial}}\hskip\parindent
The proof of Theorem \ref{thm partial} is divided into two parts. We first prove the existence of a weak solution to the FPE, and then establish the stability estimate for weak solutions.

The existence proof relies on the classical DiPerna--Lions theory, which is particularly well suited to the case of degenerate noise. In our setting, an additional difficulty arises, namely ensuring that the solution of the FPE remains a probability measure. When adapting the compactness argument, we must therefore show that no mass is lost in the limit. We first recall the compactness criterions on which the proof rests.

\begin{theorem}[Banach--Alaoglu--Bourbaki, {\cite[Theorem~3.16]{Brezis 2011}}]\label{Banach--Alaoglu--Bourbaki}
	Let $X$ be a Banach space. The closed unit ball in $X^*$ is compact in the weak-$*$ topology.
\end{theorem}

\begin{definition}[Equi-integrability]
	Let $U\subset \Rd$ be a bounded Borel set. A family
	$\{f_i\}_{i\in I}\subset L^1(U)$ is equi-integrable if, for every
	$\eta>0$, there exists $\delta>0$ such that, for every Borel set
	$E\subset U$ with $\mathcal L^d(E)\leq\delta$, one has
	\[
	\sup_{i\in I}\int_E |f_i|\,\d x\leq \eta .
	\]
\end{definition}

\begin{theorem}[Dunford--Pettis, {\cite[Theorem~1.43]{Evans 2025}}]\label{Dunford--Pettis}
	\label{Dunford-Pettis}
	Let $U\subset \Rd$ be a bounded Borel set, and let the family
	$\{f_i\}_{i\in I}\subset L^1(U)$ be bounded. Then the following are
	equivalent:
	\begin{itemize}
		\item[(i)] $\{f_i\}_{i\in I}$ is equi-integrable;
		\item[(ii)] $\{f_i\}_{i\in I}$ is weakly relatively compact in $L^1(U)$.
	\end{itemize}
\end{theorem}

Let $\phi\in C_c^\infty(\Rd)$ be a nonnegative standard mollifier. For
$0<\e<1$, set
\begin{align}\label{mollifier}
	\phi^\e(x):=\e^{-d}\phi\left(\frac{x}{\e}\right),
	\qquad x\in\Rd .
\end{align}
We define the mollified coefficients and initial datum by
$$
b^\e:=b*\phi^\e,
\qquad
a^\e:=a*\phi^\e,
\qquad
\rho_0^\e:=\rho_0*\phi^\e.
$$
We also set
$$
(b^\sigma)^\e:=b^\sigma*\phi^\e
=
b^\e-\frac12\div(a^\e).
$$
By the integrability assumptions on $b$, the boundedness of $\sigma$, and the standard properties of mollifiers, we have for any $\epsilon\in(0,1)$,
$$
b^\e \in
L^1\bigl([0,T];C_b^\infty(\Rd;\R^d)\bigr),
\quad
a^\e \in
L^1\bigl([0,T];C_b^\infty(\Rd;\R^d\times\R^d)\bigr),
\quad
\rho_0^\e\in C_b^\infty(\Rd)\cap L^1(\R^d).
$$

We consider the regularized FPE
\begin{align*}
	\begin{cases}
		\partial_t \rho^\e +\div (b^\e \rho^\e)
		=
		\dfrac{1}{2} \sum_{i,j=1}^d
		\partial_{ij} \bigl(a^\e_{ij} \rho^\e\bigr),
		&\text{in } (0,T) \times \Rd,
		\\
		\rho^\e(0,x) = \rho_0^\e(x),
		&\text{in } \Rd.
	\end{cases}
\end{align*}
By the classical parabolic theory of linear PDEs (see, e.g., \cite{Ladyzenskaja 1968}), this equation admits a unique solution
\begin{align*}\label{reg-solution}
	\rho^\e \in C([0,T];C^2_0(\Rd)),
\end{align*}where $C^2_0(\Rd)$ is the space of twice continuously differentiable functions vanishing at infinity together with their derivatives up to order $2$.
Moreover, using
the definition of $(b^\sigma)^\e$, it
can equivalently be written in divergence form as
\begin{equation}\label{reg-FPE}
	\begin{cases}
		\partial_t \rho^\e +\div \bigl((b^\sigma)^\e \rho^\e\bigr)
		=
		\dfrac{1}{2} \sum_{i,j=1}^d
		\partial_i \bigl(a^\e_{ij}\partial_j \rho^\e\bigr),
		&\text{in } (0,T) \times \Rd,
		\\
		\rho^\e(0,x) = \rho_0^\e(x),
		&\text{in } \Rd.
	\end{cases}
\end{equation}

Having made these preparations, we are now in a position to prove the existence part of Theorem \ref{thm partial}.
\begin{proof}[Proof of Theorem \ref{thm partial}(1)
	] We divide the proof into six steps.
	
	\textbf{Step I: Non-negativity and uniform a priori estimates for the regularised solutions.}	Let $\beta:\R\to\R$ be a convex smooth function with $\beta(0)=0$, and $\epsilon\in(0,1)$. At the regularised level,
	multiplying both sides of \eqref{reg-FPE} by $\beta'(\rho^\e)$, the usual chain rule
	(see, e.g., \cite{Trevisan 2016, Lions 2008}) gives
	\begin{align}\label{identity energy}
		\frac{\d}{\d t} \int_{\Rd} \beta(\rho^\e)\,\d x
		&+ \int_{\Rd}
		\bigl(\rho^\e\beta'(\rho^\e)-\beta(\rho^\e)\bigr)
		\div[(b^\sigma)^\e]\,\d x \cr
		&\qquad +
		\int_{\Rd}
		\frac{\beta''(\rho^\e)}{2}
		(\nabla\rho^\e)^*a^\e \nabla\rho^\e\,\d x
		=0 .
	\end{align}
	
	Under 
	Assumption \ref{assump1},
	suitable choices of the convex function $\beta$ in the above identity \eqref{identity energy} give the following estimates. Since the corresponding computations are standard for the
	regularised equation and can be found in \cite{Trevisan 2016, Lions 2008}, we omit the detailed calculation here and only state the resulting estimates.
	
	\begin{enumerate}
		\item[(i)] Taking smooth approximations of
		\[
		\beta(u)=u^-,
		\qquad u\in\mathbb{R},
		\]
		we obtain the non-negativity of the regularised solutions
		\[
		\rho^\e(t,x)\geq0,
		\qquad (t,x)\in[0,T]\times\Rd,
		\]
		since $\rho_0^\e\geq0$.
		
		\item[(ii)] Taking smooth approximations of
		\[
		\beta(u)=\bigl(u-\|\rho_0\|_{L_x^\infty}\bigr)^+,
		\qquad u\in\mathbb{R},
		\]
		we obtain the uniform $L^\infty$ estimate
		\begin{equation}\label{L infty}
			\sup_{\epsilon\in(0,1)}
			\|\rho^\e\|_{L^\infty([0,T]\times\Rd)}
			\leq
			\|\rho_0\|_{L_x^\infty} \exp{\left(A(T)\right)},
		\end{equation}
		where for $t\in[0,T]$,
		$$A(t):=	
		\int_0^t
		\|(\div b^\sigma)^-(s,\cdot)\|_{L_x^\infty}\,\d s
		.$$
		\item[(iii)] Taking smooth approximations of
		\[
		\beta(u)=|u|^r,
		\qquad u\in\mathbb{R},\,1\leq r<\infty,
		\]
		and applying Gr\"{o}nwall's inequality, we derive
		\begin{equation}\label{energy estimate 1}
			\sup_{\epsilon\in(0,1)}
			\|\rho^\e\|_{L_t^\infty L_x^r}
			\leq
			\|\rho_0\|_{L_x^r}
			\exp\left[\left(1-\frac1r\right)A(T)\right],
			\qquad 1\leq r<\infty .
		\end{equation}
	\end{enumerate}

	\textbf{Step II: Mass conservation for the regularised solutions.}	
	For $R>0$,
	let $\chi_R\in C_c^\infty(\Rd)$ be a cut-off function satisfying
	$$
	0\leq\chi_R\leq1,
	\qquad
	\chi_R=1 \text{ on } B_R,
	\qquad
	\chi_R=0 \text{ outside } B_{2R},
	$$
	and
	$$
	|\nabla\chi_R|\lesssim \frac{1}{R},
	\qquad
	\|\nabla^2\chi_R\|
	\lesssim  \frac{1}{R^2}.
	$$
	Multiplying \eqref{reg-FPE} by $\chi_R$ and integrating by parts, we obtain for a.e. $t\in(0,T)$,
	$$
	\frac{\d }{\d t} \int_{\Rd}\rho^\e \chi_R \, \d x=
	\int_{\Rd} \<(b^\sigma)^\e,\nabla \chi_R\>  \rho^\e\, \d x
	+\frac{1}{2}
	\int_{\Rd} \rho^\e \div((a^\e)^*
	\nabla \chi_R) \, \d x ,$$
	which gives, for every $t\in (0,T)$,
	\begin{equation}\label{pf-Thm3.4(1)-1}
		\begin{split}
			\int_{\Rd} \rho^\e(t,\cdot)\chi_R \, \d x-\int_{\Rd} \rho^\e(0,\cdot)\chi_R \, \d x
			&=
			\int_{0}^{t}	\int_{\Rd} \<(b^\sigma)^\e,\nabla \chi_R\> \rho^\e \, \d x \, \d s\cr
			&\quad +
			\frac{1}{2}\int_{0}^{t}	\int_{\Rd} \rho^\e \div((a^\e)^*\nabla\chi_R) \, \d x \,\d s.
		\end{split}
	\end{equation}
	Since $(b^\sigma)^\e\in L^1([0,T];C_b^\infty(\Rd;\R^d))$, $a^\e\in L^1([0,T];C_b^\infty(\Rd;\R^d\times\R^d))$, $\rho^\epsilon\in L^\infty([0,T];L^1(\Rd))$ with $\rho^\e\geq0$, and $\nabla \chi_R, \nabla^2 \chi_R$ are supported in $B_{2R}\setminus B_R$, the right-hand side of \eqref{pf-Thm3.4(1)-1} tends to $0$ as $R\to\infty$ by the dominated convergence theorem. Thus,
	$$\int_{\Rd} \rho^\e(t,\cdot) \, \d x=\int_{\Rd} \rho^\e(0,\cdot) \, \d x=1,$$
	so that $\rho^\e(t,\cdot)$ is a probability density for every $t\in[0,T]$.

	\textbf{Step III: Compactness and extraction of a subsequence.}	
	From \eqref{L infty} and \eqref{energy estimate 1}, $(\rho^\epsilon)_{\epsilon\in(0,1)}$ is uniformly bounded in $L^\infty([0,T]; L^r(\R^d))$ for every $r\in(1,\infty]$.
	By the Banach--Alaoglu--Bourbaki theorem (Theorem
	\ref{Banach--Alaoglu--Bourbaki}), up to extracting a subsequence, not relabeled, for every $r\in(1,\infty]$, as $\epsilon\rightarrow0^+$, we have the weak-$*$ convergence
	\begin{align}\label{w-star-limit}
		\rho^\e \xrightharpoonup{*}  \rho
		\quad\text{in } L^\infty([0, T]; L^r(\R^d)).
	\end{align}
	
	The case $r=1$ requires a separate argument,  one cannot argue directly by
	duality,
	since the dual space of $L^\infty$
	is strictly larger than $L^1$. Instead, we show that
	$(\rho^\e)_{\epsilon\in(0,1)}$ is weakly relatively compact in
	$L^1_{\loc}([0,T]\times\Rd)$. To this end, let $(\rho_0^n)_{n\geq1}\subset C_c^\infty(\Rd)$ be a nonnegative sequence converging
	to $\rho_0$ in $L^1(\Rd)$ as $n\rightarrow\infty$, and for each $n\geq1$, denote by $\rho^{n,\e}$ the
	corresponding approximating solutions constructed as above with initial datum
	$
	\rho_0^{n,\e}:=\rho_0^n*\phi^\e.
	$
	For fixed $n\geq1$ and every $1\leq r\leq \infty$, \eqref{L infty} and \eqref{energy estimate 1} give a uniform bound
	\[
	\|\rho^{n,\e}\|_{L_t^\infty L_x^r}\leq C(n,r,T),
	\]
	with $C(n,r,T)$ being a positive constant, depending on $n,r,T$ and independent of $\e$. Hence, choosing any fixed $r>1$, for every
	bounded cylinder $Q\subset [0,T]\times\Rd$ and every Borel set $E\subset Q$,
	\[
	\sup_{0<\e<1}
	\int_E \rho^{n,\e}\,\d x\,\d t
	\leq C(n,r,T) \mathcal{L}^{d+1}(E)^{1/r'},
	\]
	where $r'=\frac{r}{r-1}$.	
	Thus, for each fixed $n$, the family $\{\rho^{n,\e}\}_{0<\e<1}$ is
	equi-integrable on every bounded cylinder $Q\subset[0,T]\times\Rd$.
	
	On the other hand, by linearity, $\rho^\e-\rho^{n,\e}$ solves the same
	regularized equation with initial datum $\rho_0^\e-\rho_0^{n,\e}$. Hence, the
	corresponding $L^1$ estimate \eqref{energy estimate 1} with $r=1$ gives
	\begin{equation}\label{approx}
		\|\rho^\e-\rho^{n,\e}\|_{L_t^\infty L_x^1}
		\leq
		C_0\|\rho_0^\e-\rho_0^{n,\e}\|_{L_x^1}
		\leq
		C_0\|\rho_0-\rho_0^n\|_{L_x^1},
	\end{equation}
	where
	$
	C_0:=
	\exp(
	A(T)
	).
	$
	
	We now show that $\{\rho^\e\}_{0<\e<1}$ is equi-integrable on a fixed bounded
	cylinder $Q\subset [0,T]\times\Rd$. Let $\eta>0$. Choose $n$ so large that
	\[
	T C_0\|\rho_0-\rho_0^n\|_{L_x^1}\leq\frac{\eta}{2}.
	\]
	For this fixed $n$, choose $\delta>0$ such that, whenever $E\subset Q$ and
	$
	\mathcal{L}^{d+1}(E)
	\leq\delta$,
	\[
	\sup_{0<\e<1}
	\int_E \rho^{n,\e}\,\d x\,\d t
	\leq \frac{\eta}{2}.
	\]
	Then, combining this with
	\eqref{approx},
	\[
	\sup_{0<\e<1}
	\int_E \rho^\e\,\d x\,\d t
	\leq
	T C_0\|\rho_0-\rho_0^n\|_{L_x^1}
	+
	\sup_{0<\e<1}
	\int_E \rho^{n,\e}\,\d x\,\d t
	\leq\eta .
	\]
	Therefore, $\{\rho^\e\}_{0<\e<1}$ is equi-integrable on every bounded cylinder.
	By the Dunford--Pettis theorem (Theorem \ref{Dunford--Pettis}) and a diagonal
	extraction over an increasing sequence of bounded cylinders, we may assume (passing to a further subsequence, still denoted by $\rho^\epsilon$) that the following weak convergence holds:
	\begin{align}\label{w-limit}
		\rho^\e \rightharpoonup \rho
		\quad\text{
			in }L^1_{\loc}([0,T]\times\Rd).
	\end{align}
	The limit in \eqref{w-limit} coincides with the one in \eqref{w-star-limit} by uniqueness of weak limits. In particular,
	$$
	\rho\geq 0
	\quad\text{for a.e. }(t,x)\in[0,T]\times\Rd.
	$$
	
	\textbf{Step IV: Uniform first moment estimate.}
	We aim to establish the following estimate
	\begin{equation}\label{eq:uniform-moment}
		\sup_{0<\e<1}	\operatorname*{ess\,sup}_{t\in[0,T]}
		\int_{\Rd}|x|\rho^\e(t,x)\,\d x=:C<\infty .
	\end{equation}
	Let
	$
	\sigma^\e:=(a^\e)^{1/2}
	$
	be the symmetric non-negative square root of $a^\e$. Since
	$
	a^\e$ belongs to
	$L^1([0,T];C_b^2(\Rd;\R^{d\times d}))$
	and $a^\e$ is symmetric and non-negative definite,
	\cite[Lemma~3.2.3]{StroockVaradhan 2006} implies that $\sigma^\e$ is bounded and
	Lipschitz in $x\in\Rd$, with Lipschitz constant integrable in $t$. Hence the SDE
	$$
	\d X_t^\e
	=
	b^\e(t,X_t^\e)\,\d t
	+
	\sigma^\e(t,X_t^\e)\,\d W_t,
	\qquad
	\mathcal L_{X_0^\e}(\d x)=\rho_0^\e(x)\,\d x,
	$$
	is well-posed. By It\^o's formula, the law of $X_t^\e$ is a
	weakly
	continuous
	solution to the regularized FPE. Combined with the uniqueness of
	weakly
	continuous solutions to the FPE, see Appendix A.4 in \cite{Trevisan 2016}, this
	gives
	$$
	\mathcal L_{X_t^\e}(\d x)=\rho^\e(t,x)\,\d x
	\quad\text{for a.e. }t\in[0,T].
	$$
	Applying It\^{o}'s formula and taking expectations, we have
	$$
	\begin{aligned}
		\operatorname*{ess\,sup}_{t\in[0,T]}
		\int_{\Rd}|x|\rho^\e(t,x)\,\d x
		&=
		\sup_{t\in[0,T]}\E |X_t^\e|
		\\
		&\leq
		\int_{\Rd}|x|\rho_0^\e(x)\,\d x
		+
		\left(
		\int_0^T\int_{\Rd}
		\|\sigma^\e(s,x)\|^2
		\rho^\e(s,x)\,\d x\,\d s
		\right)^{1/2}
		\\
		&\quad
		+
		\int_0^T\int_{\Rd}
		|b^\e(s,x)|\rho^\e(s,x)\,\d x\,\d s .
	\end{aligned}
	$$
	Since $\rho_0^\e=\rho_0*\phi^\e$ and $0<\e<1$,
	$$
	\int_{\Rd}|x|\rho_0^\e(x)\,\d x
	\leq
	\int_{\Rd}|x|\rho_0(x)\,\d x+1.
	$$
	Moreover, from $(\mathcal{H}^{\sigma})$,
	$$
	\int_0^T\int_{\Rd}
	\|\sigma^\e(s,x)\|^2
	\rho^\e(s,x)\,\d x\,\d s
	\leq
	\|\sigma\|_{L_t^2L_x^\infty}^2 .
	$$
	For the drift term, using the structure
	$
	b^\e=(b_1*_1\phi^\e,b_2*\phi^\e),
	$
	mass conservation, and \eqref{L infty}, we obtain
	$$
	\begin{aligned}
		\int_0^T\int_{\Rd}
		|b^\e(s,x)|\rho^\e(s,x)\,\d x\,\d s
		&\leq
		\|b_1\|_{L_t^1L_{x_1}^\infty}
		+
		\|\rho^\e\|_{L^\infty([0,T]\times\Rd)}
		\|b_2\|_{L_t^1L_x^1}
		\\
		&\leq
		\|b_1\|_{L_t^1L_{x_1}^\infty}
		+
		C_0\|\rho_0\|_{L_x^\infty}
		\|b_2\|_{L_t^1L_x^1}.
	\end{aligned}
	$$
	Combining the above estimates proves \eqref{eq:uniform-moment}. Here the symbol $\ast_1$ denotes convolution with respect to the $x_1$-variable only. We assume that  the mollifier splits as $\phi(x_1, x_2) = \phi_1(x_1)\phi_2(x_2)$ for $(x_1,x_2)\in\R^{d_1}\times\R^{d_2}$, and set $\phi_1^\epsilon(x_1) = \epsilon^{-d_1}\phi_1(x_1/\epsilon)$, where $\phi_1,\phi_2$ are nonnegative standard mollifiers on $\R^{d_1}$ and $\R^{d_2}$, respectively. Then $b_1\ast_1\phi^\epsilon:=b_1\ast\phi_1^\epsilon$.
	
	\textbf{Step V: 	
		Recovery of the total mass of the limit density.} For every
	nonnegative $\psi\in L^\infty([0,T])$ and the cut-off function $\chi_R$, the local
	weak convergence in \eqref{w-limit} gives
	\begin{align*}
		\int_0^T\int_{\Rd}
		\rho(t,x)\psi(t)\chi_R(x)\,\d x\,\d t
		&=
		\lim_{\e\to0}
		\int_0^T\int_{\Rd}
		\rho^\e(t,x)\psi(t)\chi_R(x)\,\d x\,\d t\cr
		&\leq
		\int_0^T\psi(t)\,\d t .
	\end{align*}
	Letting $R\to\infty$ yields
	$$
	\int_{\Rd}\rho(t,x)\,\d x\leq1,
	\quad\text{for a.e. }t\in[0,T].
	$$
	Thus, for a.e. $t\in[0,T]$, $\rho(t,\cdot)$
	is a sub-probability density.
	
	Conversely, by \eqref{eq:uniform-moment}, for a.e. $t\in[0,T]$,
	$$
	\int_{\Rd}
	\rho^\e(t,x)(1-\chi_R(x))\,\d x
	\leq
	\int_{|x|\geq R}
	\rho^\e(t,x)\,\d x
	\leq
	\frac1R
	\int_{\Rd}
	|x|\rho^\e(t,x)\,\d x
	\leq
	\frac{C}{R}.
	$$
	Since $\rho^\e(t,\cdot)$ is a probability density,
	$$
	\int_{\Rd}
	\rho^\e(t,x)\chi_R(x)\,\d x
	=
	1-
	\int_{\Rd}
	\rho^\e(t,x)(1-\chi_R(x))\,\d x
	\geq
	1-\frac{C}{R}.
	$$
	Multiplying by $\psi(t)$, integrating over $(0,T)$, and passing to the limit
	$\e\to0$, we obtain
	$$
	\int_0^T\int_{\Rd}
	\rho(t,x)\psi(t)\chi_R(x)\,\d x\,\d t
	\geq
	\left(1-\frac{C}{R}\right)
	\int_0^T\psi(t)\,\d t .
	$$
	Letting $R\to\infty$ gives
	$$
	\int_0^T\int_{\Rd}
	\rho(t,x)\psi(t)\,\d x\,\d t
	\geq
	\int_0^T\psi(t)\,\d t .
	$$
	Together with the sub-probability estimate, this proves
	$$
	\int_{\Rd}\rho(t,x)\,\d x=1,
	\quad\text{for a.e. }t\in[0,T].
	$$
	Consequently,
	$$
	\rho\in L^\infty([0,T]; L^1(\R^d)\cap L^\infty(\R^d)),
	$$
	and $\rho(t,\cdot)$ is a probability density for a.e. $t\in[0,T]$.
	
	\textbf{Step VI: Verification that the limit 
		$\rho$ satisfies \eqref{weak solution integral}.} Let
	$\varphi\in C_c^{1,2}([0,T)\times\Rd)$. The regularized solutions satisfy
	$$
	\int_0^T\int_{\Rd}
	\rho^\e
	\left[
	\partial_t\varphi
	+
	b^\e\cdot\nabla\varphi
	+
	\frac12 \tr\bigl(a^\e \nabla^2\varphi\bigr)
	\right]
	\,\d x\,\d t
	+
	\int_{\Rd}\rho_0^\e(x)\varphi(0,x)\,\d x
	=0 .
	$$
	Using the local convergence of the mollified coefficients, the uniform
	$L^\infty$ bound on $\rho^\e$, the weak-$*$ convergence of $\rho^\e$ in
	$L^\infty([0,T]\times \R^d)$, and the weak convergence in
	$L^1_{\loc}([0,T]\times\R^d)$, we may pass to the limit
	in each term. Since $\rho_0^\e\to\rho_0$ in $L^1(\Rd)$ as $\epsilon\rightarrow0^+$, we obtain
	$$
	\int_0^T\int_{\Rd}
	\rho
	\left[
	\partial_t\varphi
	+
	b\cdot\nabla\varphi
	+
	\frac12 \tr\bigl(a \nabla^2\varphi\bigr)
	\right]
	\,\d x\,\d t
	+
	\int_{\Rd}\rho_0(x)\varphi(0,x)\,\d x
	=0 .
	$$
	By Remark \ref{narrowly continuous}, the measure-valued map $\rho(t,x)\,\d x$ admits a unique weakly continuous representative with initial value $\rho_0(x)\,\d x$. This completes the construction of a weak solution to the FPE.
\end{proof}

\begin{remark}
	By Theorem \ref{thm partial}(1) and the superposition principle Theorem \ref{sp}, the FPE solution induces a weak solution to the SDE \eqref{SDE}, which also yields a construction of weak solutions to the associated degenerate SDE.
\end{remark}

We now turn to the second part of Theorem \ref{thm partial}, the quantitative stability
estimate for weak solutions. The proof relies on several auxiliary lemmas, which
we state first.

Since our analysis concerns stability under rough coefficients, we use a standard
mollification procedure. In particular, we shall use a lemma due to Crippa and
Ligabue \cite{Crippa 2020}. Under fractional Sobolev regularity, it provides
quantitative estimates both for the convergence of a mollified function to the
original one and for the blow-up rate of its derivatives.

\begin{lemma} \label{mollified function rate}
	Let $s\in(0,1)$ and $p\in[1,\infty)$. Then, there exists a constant $C>0$ such that for any $\epsilon\in(0,1)$ and $f\in W^{s,p}(\R^d)$,
	$$
	\|f- f*\phi^\e\|_{L^p(\R^d)}
	\leq
	C\e^s\|f\|_{W^{s,p}(\R^d)},
	\qquad
	\|\nabla (f*\phi^\e)\|_{L^p(\R^d)}
	\leq
	C\e^{s-1}\|f\|_{W^{s,p}(\R^d)},
	$$
	where $\phi^\epsilon$ is given by \eqref{mollifier}.
\end{lemma}

To control the error terms below, we shall use the Hardy--Littlewood maximal
function. For $f\in L_{\loc}^1(\Rd)$, set
$$
\M f(x):=
\sup_{r>0}
\frac{1}{\mathcal L^d(B_r(x))}
\int_{B_r(x)} |f(z)|\,\d z,
\quad x\in\Rd.
$$
It is well known that $\M$ is of weak type $(1,1)$, i.e., there exists a constant $c_d>0$ depending only on $d$ such that
\begin{equation}\label{weak (1,1)}
	\|\M f\|_{L^{1,\infty}(\Rd)}
	\leq c_d\|f\|_{L^1(\R^d)},\quad f\in L^1(\R^d).
\end{equation}

Moreover, for $1<r\leq\infty$,  $\M$  satisfies the strong $(r,r)$ estimate, i.e., there exists a constant $c_{d,r}>0$ depending only on $d,r$ such that
$$
\|\M f\|_{L^r(\R^d)}
\leq c_{d,r}\|f\|_{L^r(\R^d)},\quad f\in L^r(\R^d).
$$

We next recall a Morrey-type inequality, which bounds the increments of a Sobolev
function in terms of the maximal function of its gradient; see
\cite[Appendix]{Fang 2010} or the more precise version \cite[Lemma 3.1]{Jabin2010} for a proof.
\begin{lemma}\label{Morrey-type estimate}
	There exists a constant $C_d> 0$ such that, for any $f\in W_{\loc}^{1,1}(\Rd)$,
	$$
	|f(x)-f(y)|
	\leq C_d|x-y|
	\left(
	\M(\nabla f)(x)+\M(\nabla f)(y)
	\right),
	$$
	for a.e. $x,y\in\Rd$.
\end{lemma}
The following lemmas give a variant of this Morrey-type estimate, together with
the corresponding $L^r$-estimate for the modified maximal function. These results
were originally proved in \cite{Crippa 2020}; for completeness, we include a brief
proof.

\begin{lemma}[A variant maximal function estimate]\label{A variant maximal function estimate}
	There exists a constant $C_d>0$ such that, for every $f\in W_{\loc}^{1,1}(\Rd)$ and for a.e. $x,y\in\Rd$,
	$$
	|f(x)-f(y)|
	\leq
	C_d |B^{-1}(x-y)|
	\left(
	\widetilde{\M}(\nabla f)(x)
	+
	\widetilde{\M}(\nabla f)(y)
	\right),
	$$
	where the matrix $B$ is given by \eqref{diag-B}, and
	$$
	\widetilde{\M}(\nabla f)(x)
	:=
	\M\left(
	\sum_{i=1}^{d}
	|\partial_i f(B\cdot)|B_{ii}
	\right)(B^{-1}x).
	$$
\end{lemma}

\begin{proof}
	Set $\widetilde f(x):=f(Bx)$. By Lemma \ref{Morrey-type estimate}, for a.e. $\xi,\eta\in\Rd$,
	$$
	|\widetilde f(\xi)-\widetilde f(\eta)|
	\leq
	C_d|\xi-\eta|
	\left(
	\M(\nabla \widetilde f)(\xi)
	+
	\M(\nabla \widetilde f)(\eta)
	\right).
	$$
	Since
	$$
	\M(\nabla \widetilde f)(\xi)
	\leq
	\M\left(
	\sum_{i=1}^{d}
	|\partial_i f(B\cdot)|B_{ii}
	\right)(\xi),
	$$
	we obtain, for a.e. $\xi,\eta\in\Rd$,
	$$
	|f(B\xi)-f(B\eta)|\leq
	C_d|\xi-\eta|\left[\M\left(\sum_{i=1}^{d}|\partial_i f(B\cdot)|B_{ii}\right)(\xi)
	+
	\M\left(\sum_{i=1}^{d}|\partial_i f(B\cdot)|B_{ii}\right)(\eta)\right].
	$$
	Letting $\xi=B^{-1}x$ and $\eta=B^{-1}y$ completes the proof.
\end{proof}

\begin{lemma}[Operator bounds]\label{Operator bounds}
	Let $\widetilde{\M}(\nabla f)$ be defined as in Lemma \ref{A variant maximal function estimate}. Then the following estimates hold.
	\begin{itemize}
		\item[1.] If $p=1$ and $|\nabla f|\in L^1(\R^d)$,
		then
		$$
		\|\widetilde{\M}(\nabla f)\|_{L^{1,\infty}(\Rd)}
		\lesssim
		\delta_1 \sum_{i=1}^{d_1}\|\partial_i f\|_{L^1(\R^d)}
		+
		\delta_2 \sum_{i=d_1+1}^{d}\|\partial_i f\|_{L^1(\R^d)}.
		$$
		
		\item[2.] If $1<p\leq \infty$ and 
		$|\nabla f|\in L^p(\R^d)$,
		then
		$$
		\|\widetilde{\M}(\nabla f)\|_{L^p(\R^d)}
		\lesssim
		\delta_1 \sum_{i=1}^{d_1}\|\partial_i f\|_{L^p(\R^d)}
		+
		\delta_2 \sum_{i=d_1+1}^{d}\|\partial_i f\|_{L^p(\R^d)}.
		$$
	\end{itemize}
\end{lemma}

\begin{proof}
	As in Lemma \ref{A variant maximal function estimate}, consider the rescaled
	function $\widetilde{f}(x)=f(Bx)$. Adapting the argument behind \eqref{weak (1,1)}, we obtain
	\begin{equation*}
		\begin{split}
			\|\widetilde{\M}(\nabla f)\|_{L^{1,\infty}(\Rd)}&= \left\|\M\left(\sum_{i=1}^{d}|\partial_i
			\widetilde{f}|(B^{-1}\cdot)\right)\right\|_{L^{1,\infty}(\Rd)}\\
			&\lesssim \|\sum_{i=1}^{d}|\partial_i \widetilde{f}|(B^{-1}\cdot)\|_{L^1(\R^d)}\\
			&\lesssim \sum_{i=1}^{d}\left\|\left(\partial_i f (B\cdot)B_{ii}\right)(B^{-1}\cdot)\right\|_{L^1(\R^d)}\\
			&=\sum_{i=1}^{d} B_{ii}\|\partial_i f\|_{L^1(\R^d)}.
		\end{split}
	\end{equation*}
	The proof of conclusion 2 is similar and is omitted.
\end{proof}

With these preparations, we now prove the second part of Theorem \ref{thm partial}.
\begin{proof}[Proof of Theorem \ref{thm partial}(2)]
	To circumvent the non-differentiability of
	$$
	x\mapsto \log\left(1+\sqrt{|B^{-1}x|^2}\right)
	$$
	at $x=0$, we use a standard regularization argument. For $\alpha>0$, introduce
	the smooth approximation
	$$
	x\mapsto
	\log\left(1+\sqrt{|B^{-1}x|^2+\alpha}\right),
	$$
	which is smooth on $\Rd$ and converges pointwise to
	$
	\log\left(1+\sqrt{|B^{-1}x|^2}\right)
	$
	as $\alpha\to0$.
	
	By Theorem \ref{thm partial}(1), the weak solutions $\rho$ and $\hat{\rho}$ to the FPE \eqref{FPE} with corresponding coefficients exist. By the superposition principle (Theorem 2.5), each of them gives rise to a weak solution of the SDE \eqref{SDE}, although they may be defined on different probability spaces with different Brownian motions. Applying Lemma \ref{lem-common} to these two weak solutions, we may assume, without loss of generality, that they are realized on a common filtered probability space and driven by the same Brownian motion $(W_t)_{t\geq0}$. Moreover, by the same lemma, the initial coupling can be chosen arbitrarily, and we take it to be optimal with respect to $\D_{\delta_1,\delta_2}$. With these preparations, we consider the coupled SDEs:
	$$
	\d X_t=b(t,X_t)\,\d t+\sigma(t,X_t)\,\d W_t,
	\qquad
	\L_{X_0}=\nu,
	$$
	and
	$$
	\d \hat X_t=\hat b(t,\hat X_t)\,\d t+\hat\sigma(t,\hat X_t)\,\d W_t,
	\qquad
	\L_{\hat X_0}=\hat\nu,
	$$
	with $\nu=\rho_0\d x$ and $\hat{\nu}=\hat{\rho}_0\d x$. Moreover, the superposition principle gives the time marginals
	$$
	\L_{X_t}(\d x)=\rho(t,x)\,\d x,
	\quad
	\L_{\hat X_t}(\d x)=\hat\rho(t,x)\,\d x,
	\qquad t\in(0,T).
	$$
	We also choose the initial values $(X_0,\hat X_0)$ to be optimally coupled, namely
	\begin{align}\label{pf-thm partial(2)-1}
		\E \log\left(1+\sqrt{|B^{-1}(X_0-\hat X_0)|^2}\right)
		=\D_{\delta_1,\delta_2}(\nu,\hat\nu).
	\end{align}

	Recall that
	$$
	\tau_R:=\inf\left\{t\in[0,T]:\ |X_t|\vee|\hat X_t|\geq R\right\},\quad R>0.
	$$
	By It\^o's formula applied to the above SDEs, we obtain
	\begin{align}\label{pf-thm partial(2)-2}
		&\E\log\left(
		1+\sqrt{
			|B^{-1}(X_{t\wedge\tau_R}-\hat X_{t\wedge\tau_R})|^2+\alpha
		}
		\right)
		\cr
		&\leq
		\E\log\left(
		1+\sqrt{|B^{-1}(X_0-\hat X_0)|^2+\alpha}
		\right)
		\cr
		&\qquad
		+
		\E\left(
		\int_0^{t\wedge\tau_R}
		\frac{
			|B^{-1}(b(s,X_s)-\hat b(s,\hat X_s))|
		}{
			1+|B^{-1}(X_s-\hat X_s)|
		}
		\,\d s
		\right)
		\cr
		&\qquad
		+
		\frac{1}{2}\E\left(
		\int_0^{t\wedge\tau_R}
		\frac{
			\|B^{-1}(\sigma(s,X_s)-\hat\sigma(s,\hat X_s))\|^2
		}{
			\sqrt{|B^{-1}(X_s-\hat X_s)|^2+\alpha}
			+
			|B^{-1}(X_s-\hat X_s)|^2+\alpha
		}
		\,\d s
		\right).
	\end{align}
	We set
	$$
	\mathrm{I}_1:=
	\E\left(
	\int_0^{t\wedge\tau_R}
	\frac{
		|B^{-1}(b(s,X_s)-\hat b(s,\hat X_s))|
	}{
		1+|B^{-1}(X_s-\hat X_s)|
	}
	\,\d s
	\right),
	$$
	and
	$$
	\mathrm{I}_2:=
	\E\left(
	\int_0^{t\wedge\tau_R}
	\frac{
		\|B^{-1}(\sigma(s,X_s)-\hat\sigma(s,\hat X_s))\|^2
	}{
		\sqrt{|B^{-1}(X_s-\hat X_s)|^2+\alpha}
		+
		|B^{-1}(X_s-\hat X_s)|^2+\alpha
	}
	\,\d s
	\right).
	$$
	It suffices to estimate $\mathrm{I}_1$ and $\mathrm{I}_2$.

	\textbf{Estimate for $\mathrm{I}_1$.} We split the estimate of $\mathrm{I}_1$ into four steps.
	
	\textbf{Step 1: Mollification of the drift coefficient $b$.}
	Denote by $\b_2^\e$ the regularization of $b_2$ in the $x_1$-variable,
	$$
	\b_2^\e(t,x_1,x_2)
	:=
	b_2(t,\cdot,x_2)*_1\phi^\e(x_1),
	\quad
	t\in[0,T],\ x_1\in\R^{d_1},\ x_2\in\R^{d_2},
	$$
	where $*_1$ denotes convolution in the $x_1$-variable only, with $\phi^\e$
	understood as a standard mollifier on $\R^{d_1}$, and set	
	$$
	\b^\e:=(b_1,\b_2^\e).
	$$
	
	\textbf{Step 2: Intermediate estimate via mollification.}
	Let $\Omega_{T,R}:=[0,T]\times B_R\times B_R$. Introducing the mollified
	drift as an intermediate term, we obtain
	\begin{align*}
		\mathrm{I}_1 &\leq \int_{\Omega_{T,R}} \frac{|B^{-1}(b(s,x)-\hat{b}(s,y))|}{1+|B^{-1}(x-y)|} \, \pi_s(\d x, \d y)\, \d s \\
		& \leq
		\int_{\Omega_{T,R}} \left\{\frac{|B^{-1}(b(s,x)-\b^\e(s,x))|}{1+|B^{-1}(x-y)|}+
		\frac{|B^{-1}(\b^\e(s,x)-\hat{\b}^\e(s,y))|}{1+|B^{-1}(x-y)|}             \right\} \,\pi_s(\d x, \d y) \, \d s  \\
		&\quad +  \int_{\Omega_{T,R}}\frac{|B^{-1}(\hat{\b}^\e(s,y)-\hat{b}(s,y))|}{1+|B^{-1}(x-y)|} \,\pi_s(\d x, \d y) \,\d s \\
		&\lesssim \frac{1}{\delta_2}\|b_2-\b_2^\e\|_{L^1_tL_x^1}+\frac{1}{\delta_2}\|\hat{b}_2-\hat{\b}_2^\e\|_{L^1_tL_x^1}
		+\int_{\Omega_{T,R}} \frac{|B^{-1}(\b^\e(s,x)-\hat{\b}^\e(s,y))|}{1+|B^{-1}(x-y)|} \,\pi_s(\d x, \d y) \,\d s,
	\end{align*}
	where $\pi_s=\L_{(X_s,\hat X_s)}$ stands for the joint distribution of $(X_s,\hat X_s)$ for any $s\in[0,T]$. The last inequality uses the boundedness of the densities $\rho$ and $\hat\rho$.
	
	Set
	$$\mathrm{I}_{1,1}:= \int_{\Omega_{T,R}}  \frac{|B^{-1}(\b^\e(s,x)-\hat{\b}^\e(s,y))|}{1+|B^{-1}(x-y)|} \, \pi_s(\d x, \d y) \, \d s.$$
	Then, for any $0<\delta_1<\delta_2$,
	\begin{equation*}
		\begin{aligned}
			\mathrm{I}_{1,1}
			&\leq
			\int_{\Omega_{T,R}}
			\left\{
			\frac{|B^{-1}(\b^\e(s,y)-\hat{\b}^\e(s,y))|}{1+|B^{-1}(x-y)|}
			+\frac{|B^{-1}(\b^\e(s,x)-{\b}^\e(s,y))|}{1+|B^{-1}(x-y)|}
			\right\}
			\,\pi_s(\d x,\d y)\,\d s
			\\
			&\leq
			\frac{1}{\delta_1}
			\int_{\Omega_{T,R}}
			|b_1(s,y_1)-\hat b_1(s,y_1)|
			\hat\rho_s(y_1,y_2)
			\,\d y_1\,\d y_2\,\d s
			+
			\frac{1}{\delta_2}
			\|b_2-\hat b_2\|_{L_t^1L_x^1}
			\\
			&\quad
			+
			\int_{\Omega_{T,R}}
			\frac{|B^{-1}(\b^\e(s,x)-{\b}^\e(s,y))|}{1+|B^{-1}(x-y)|}
			\,\pi_s(\d x,\d y)\,\d s
			\\
			&\lesssim
			\int_{\Omega_{T,R}}
			\frac{|B^{-1}(\b^\e(s,x)-{\b}^\e(s,y))|}{1+|B^{-1}(x-y)|}
			\,\pi_s(\d x,\d y)\,\d s
			+
			\frac{1}{\delta_1}
			\|b_1-\hat b_1\|_{L^1([0,T];L^p(B_R))}
			\\
			&\quad
			+
			\frac{1}{\delta_1}
			\|b_2-\hat b_2\|_{L_t^1L_x^1},
		\end{aligned}
	\end{equation*}
	where in the last inequality we used $0<1/\delta_2<1/\delta_1$ and the
	boundedness of the densities $\rho$ and $\hat\rho$. It remains to estimate
	$$\mathrm{I}_{1,1,1}:=\int_{\Omega_{T,R}} \frac{|B^{-1}(\b^\e(s,x)-{\b}^\e(s,y))|}{1+|B^{-1}(x-y)|}\,\pi_s(\d x, \d y) \, \d s.$$
	Applying Lemma \ref{A variant maximal function estimate} to the components
	$b_1$ and $\b^\e_2$, together with the subadditivity of $\widetilde{\M}$, we
	obtain, for a.e. $x,y \in \Rd$,
	$$
	\frac{|b_1(t,x)-b_1(t,y)|}{|B^{-1}(x-y)|}
	\lesssim
	\widetilde{\M}(\nabla_{x_1}b_1)(t,x)
	+
	\widetilde{\M}(\nabla_{x_1}b_1)(t,y),
	$$
	and
	\begin{equation*}
		\begin{aligned}
			\frac{|\b_2^\e(t,x)-\b_2^\e(t,y)|}{|B^{-1}(x-y)|}
			&\lesssim
			\widetilde{\M}(\nabla\b_2^\e)(t,x)
			+
			\widetilde{\M}(\nabla\b_2^\e)(t,y)
			\\
			&\lesssim
			\widetilde{\M}(\nabla_{x_1}\b_2^\e)(t,x)
			+
			\widetilde{\M}(\nabla_{x_1}\b_2^\e)(t,y)
			\\
			&\quad
			+
			\widetilde{\M}(\nabla_{x_2}\b_2^\e)(t,x)
			+
			\widetilde{\M}(\nabla_{x_2}\b_2^\e)(t,y).
		\end{aligned}
	\end{equation*}
	Since the marginals of $\pi_t$ are absolutely continuous with respect to the
	Lebesgue measure, these pointwise estimates hold for a.e. $t\in[0,T]$ and
	$\pi_t$-a.e. $(x,y)\in\Rd\times\Rd$. Therefore
	\begin{equation*}
		\mathrm{I}_{1,1,1}\leq \int_{\Omega_{T,R}}\psi_1 \, \pi_s(\d x, \d y)\, \d s+\int_{\Omega_{T,R}}\psi_2 \, \pi_s(\d x, \d y) \, \d s+\int_{\Omega_{T,R}}\psi_3\, \pi_s(\d x, \d y) \,\d s,
	\end{equation*}
	where
	$$
	\psi_1(s,x,y):=
	\min \left\{
	|B^{-1}(\b^\e(s,x)-{\b}^\e(s,y))|,
	\frac{1}{\delta_1}
	\left(
	\widetilde{\M}(\nabla_{x_1}b_1)(s,x)
	+
	\widetilde{\M}(\nabla_{x_1}b_1)(s,y)
	\right)
	\right\},
	$$
	$$
	\psi_2(s,x,y):=
	\min \left\{
	|B^{-1}(\b^\e(s,x)-{\b}^\e(s,y))|,
	\frac{1}{\delta_2}
	\left(
	\widetilde{\M}(\nabla_{x_1}\b_2^\e)(s,x)
	+
	\widetilde{\M}(\nabla_{x_1}\b_2^\e)(s,y)
	\right)
	\right\},
	$$
	and
	$$
	\psi_3(s,x,y):=
	\min \left\{
	|B^{-1}(\b^\e(s,x)-{\b}^\e(s,y))|,
	\frac{1}{\delta_2}
	\left(
	\widetilde{\M}(\nabla_{x_2}\b_2^\e)(s,x)
	+
	\widetilde{\M}(\nabla_{x_2}\b_2^\e)(s,y)
	\right)
	\right\}.
	$$
	
	\textbf{Step 3: An estimate for the maximal function.}
	Define a finite measure $\mu$ on $\Omega_T:=[0,T]\times\Rd \times \Rd$ by
	$$\mu(\d t, \d x, \d y):=\mathbbm{1}_{[0,T]}(t)\mathbbm{1}_{B_R\times B_R}(x,y) \pi_t(\d x, \d y) \d t.$$
	Using the second term in the minimum and the operator bounds in Lemma
	\ref{Operator bounds}, since $p\in(1,\infty]$, we deduce
	\begin{equation*}
		\begin{aligned}
			\|\psi_1\|_{L^1(\Omega_T,\mu)}
			&\lesssim
			\frac{1}{\delta_1}
			\|\widetilde{\M}(\nabla_{x_1}b_1)\|_{L^1(\Omega_T,\mu)}
			\\
			&\lesssim
			\frac{1}{\delta_1}
			\int_0^T
			\int_{B_R^{d_1}}
			\widetilde{\M}(\nabla_{x_1}b_1)(t,x_1)
			\,\d x_1\,\d t
			\\
			&\lesssim
			\frac{1}{\delta_1}
			\int_0^T
			\left(
			\int_{B_R^{d_1}}
			\widetilde{\M}(\nabla_{x_1}b_1)(t,x_1)^p
			\,\d x_1
			\right)^{\frac{1}{p}}
			\,\d t
			\\
			&\lesssim
			\frac{1}{\delta_1}
			\|\widetilde{\M}(\nabla_{x_1}b_1)\|_{L_t^1L_{x_1}^p}
			\\
			&\lesssim
			\|\nabla_{x_1}b_1\|_{L_t^1L_{x_1}^p}.
		\end{aligned}
	\end{equation*}
	Proceeding similarly,
	\begin{equation*}
		\begin{split}
			\|\psi_3\|_{L^1(\Omega_T,\mu)} & \lesssim \|\nabla_{x_2}b_2\|_{L^1([0,T]\times \R^{d_1}; L_{x_2}^p)}.\\
		\end{split}
	\end{equation*}
	
	The estimate of $\psi_2$ requires more care, due to the lack of regularity of
	$b_2$ in the $x_1$-variable. On the one hand, Lemma \ref{mollified function rate}, the weak $(1,1)$ bound for $\widetilde{\M}$, and the uniform boundedness of the densities give
	\begin{equation*}
		\begin{split}
			\|\psi_2\|_{L^{1,\infty}(\Omega_T,\mu)}
			&\leq
			\frac{2}{\delta_2}\left(\|\rho\|_{L_t^\infty L_x^\infty}+\|\hat{\rho}\|_{L_t^\infty L_x^\infty}\right) \int_{0}^{T} \|\widetilde{\M}(\nabla_{x_1}\b_2^\e)(t,\cdot)\|_{L^{1,\infty}(\Rd)} \, \d t\\
			& \lesssim
			\frac{\delta_1}{\delta_2}\|\nabla_{x_1}\b_2^\e\|_{L_t^1L_x^1}\\
			& \lesssim
			\frac{\delta_1}{\delta_2} \e^{s-1} \|b_2\|_{L^1([0,T]\times \R^{d_2}; W_{x_1}^{s,1})}\\
			& \lesssim
			\frac{\delta_1}{\delta_2} \e^{s-1}.
		\end{split}
	\end{equation*}
	The right-hand side diverges as $\e \to 0^+$ whenever $0<s<1$, which highlights
	the importance of a careful choice of $\delta_1$ and $\delta_2$. On the other
	hand, by the local $L^p$-integrability of $b$ on $[0,T]\times \Rd$,
	\begin{equation*}
		\begin{split}
			\|\psi_2\|_{L^p(\Omega_T,\mu)} &  \lesssim \frac{1}{\delta_1}\|\b^\e\|_{L^p([0,T]\times B_R)}\\
			&  \lesssim \frac{1}{\delta_1}\|b\|_{L^p([0,T]\times B_{R+1})}\\
			&  \lesssim \frac{1}{\delta_1}.
		\end{split}
	\end{equation*}
	Applying the interpolation inequality \eqref{interp 2} and the bound
	$\|\psi_2\|_{L^{p,\infty}(\Omega_T,\mu)}\leq\|\psi_2\|_{L^{p}(\Omega_T,\mu)}$, we obtain
	$$\|\psi_2\|_{L^1(\Omega_T,\mu)}\leq \frac{p}{p-1}\|\psi_2\|_{L^{1,\infty}(\Omega_T,\mu)}\left[1+\log\left(\frac{\mu(\Omega_T)^{1-1/p}\|\psi_2\|_{L^p(\Omega_T,\mu)}}{\|\psi_2\|_{L^{1,\infty}(\Omega_T,\mu)}}\right)\right].$$
	Since
	$$
	\|\psi_2\|_{L^{1,\infty}(\Omega_T,\mu)}
	\leq
	\mu(\Omega_T)^{1-\frac{1}{p}}
	\|\psi_2\|_{L^{p,\infty}(\Omega_T,\mu)}
	\leq
	\mu(\Omega_T)^{1-\frac{1}{p}}
	\|\psi_2\|_{L^p(\Omega_T,\mu)},
	$$
	the right-hand side is nondecreasing as a function of
	$\|\psi_2\|_{L^{1,\infty}(\Omega_T,\mu)}$. Thus,
	\begin{equation*}
		\|\psi_2\|_{L^1(\Omega_T,\mu)}   \lesssim \frac{\delta_1}{\delta_2} \e^{s-1} \left[1+\log\left(\frac{1}{\delta_1 \frac{\delta_1}{\delta_2}\e^{s-1}}\right) \right].
	\end{equation*}
	
	Combining all the above estimates, we derive
	\begin{align*}
		\mathrm{I}_1 &\lesssim \frac{1}{\delta_2}\|b_2-\b_2^\e\|_{L^1_tL_x^1}+\frac{1}{\delta_2}\|\hat{b}_2-\hat{\b}_2^\e\|_{L^1_tL_x^1}+\frac{1}{\delta_1} \|b_1-\hat{b}_1\|_{L^1([0,T];L^p(B_R))}\\
		&\quad +\frac{1}{\delta_1}\|b_2-\hat{b}_2\|_{L_t^1L_x^1}  +  \frac{\delta_1}{\delta_2} \e^{s-1} \left[1+\log\left(\frac{1}{\delta_1 \frac{\delta_1}{\delta_2}\e^{s-1}}\right) \right]+1\\
		& \lesssim \frac{1}{\delta_1} \|b_1-\hat{b}_1\|_{L^1([0,T];L^p(B_R))}+\frac{1}{\delta_1}\|b_2-\hat{b}_2\|_{L_t^1L_x^1}+\frac{\e^s}{\delta_2}\\
		&\quad +\frac{\delta_1}{\delta_2} \e^{s-1} \left[1+\log\left(\frac{\delta_2}{\delta_1^2 \e^{s-1}}\right) \right]+1.
	\end{align*}

	\textbf{Step 4: Choice of the parameters.}
	Set
	$$\beta:=\frac{\delta_1}{\delta_2} \ll 1,$$
	and let $\e=\beta$. Then
	$$\frac{\e^s}{\delta_2}+\frac{\delta_1}{\delta_2} \e^{s-1} \left[1+\log\left(\frac{\delta_2}{\delta_1^2  \e^{s-1}}\right) \right]=\frac{\beta^s}{\delta_2}+\beta^s\left[1+\log \left(\frac{1}{\delta_2 \beta^{s+1}}\right)\right].$$
	For $\delta_2\in (0,1)$,
	\begin{equation*}
		\begin{aligned}
			\frac{\beta^s}{\delta_2}
			+
			\beta^s
			\left[
			1+\log \left(\frac{1}{\delta_2 \beta^{s+1}}\right)
			\right]
			&=
			\beta^s
			\left[
			\frac{1}{\delta_2}
			+
			1
			-
			\log \delta_2
			-
			(s+1)\log \beta
			\right]
			\\
			& \leq
			\beta^s
			\left[
			\frac{1}{\delta_2}
			+
			1
			-
			\log \delta_2
			\right]
			+
			2\beta^s(-\log \beta).
		\end{aligned}
	\end{equation*}
	Since
	$$\lim_{\beta \to 0}2\beta^s(-\log\beta)=0,$$
	there exists a constant $\beta_0\in (0,1)$ such that, for $\beta\in (0, \beta_0]$,
	$$2\beta^s(-\log\beta)<1.$$
	Consequently, there exists
	$$\delta_1\leq \min\left\{\delta_2\beta_0, \delta_2 \left(\frac{1}{\left[\frac{1}{\delta_2}+1-\log(\delta_2)\right]^{1/s}}\right)\right\}<1$$
	such that
	\begin{align}\label{pf-thm partial(2)-3}
		\mathrm{I}_1\lesssim \frac{1}{\delta_1} \|b_1-\hat b_1\|_{L^1([0,T]; L^p(B_R))}+\frac{1}{\delta_1}\|b_2-\hat{b}_2\|_{L_t^1L_x^1}+1.
	\end{align}

	\textbf{Estimate for $\mathrm{I}_2$.}
	Since the diffusion coefficient is the Sobolev-regular part, we estimate
	$\mathrm{I}_2$ directly:
	\begin{align*}
		\mathrm{I}_2
		&\lesssim
		\int_{\Omega_{T,R}}
		\frac{
			\|B^{-1}(\sigma(s,x)-\hat\sigma(s,y))\|^2
		}{
			\sqrt{|B^{-1}(x-y)|^2+\alpha}
			+
			|B^{-1}(x-y)|^2+\alpha
		}
		\,\pi_s(\d x,\d y)\,\d s
		\\
		&\lesssim
		\int_{\Omega_{T,R}}
		\frac{
			\|B^{-1}(\sigma(s,x)-\sigma(s,y))\|^2
		}{
			|B^{-1}(x-y)|^2
		}
		\,\pi_s(\d x,\d y)\,\d s
		\\
		&\quad
		+
		\int_{\Omega_{T,R}}
		\frac{
			\|B^{-1}(\sigma(s,y)-\hat\sigma(s,y))\|^2
		}{
			\sqrt{\alpha}
		}
		\,\pi_s(\d x,\d y)\,\d s
		\\
		&\lesssim \frac{1}{\delta_1^2\sqrt{\alpha}}\int_0^T\int_{B_R}\|\sigma(s,y)-\hat\sigma(s,y)\|^2\hat{\rho}(y)\,\d y\,\d t\\
		&\quad + \int_{\Omega_{T,R}}\frac{\|B^{-1}(\sigma(s,x)-\sigma(s,y))\|^2}{|B^{-1}(x-y)|^2}\,\pi_s(\d x,\d y)\,\d s\\
		&\lesssim
		\frac{1}{\delta_1^2\sqrt{\alpha}}
		\|\sigma-\hat\sigma\|^2_{L^2([0,T];L^{2r}(B_R))}
		\\
		&\quad
		+
		\int_{\Omega_{T,R}}
		\frac{
			\|B^{-1}(\sigma(s,x)-\sigma(s,y))\|^2
		}{
			|B^{-1}(x-y)|^2
		}
		\,\pi_s(\d x,\d y)\,\d s .
	\end{align*}
	Here, for the penultimate inequality,  we use that $\rho\in L_t^\infty L_x^1$ and $\delta_1\in(0,\delta_2)$,  and for the last inequality, we applied H\"{o}lder's inequality since $\sigma,\hat{\sigma} \in L^2_tL^\infty_x$  and $\hat{\rho}\in L_t^\infty L_x^{r'}$ with $r'$ being the conjugate exponent of $r$.

	Set
	$$
	\mathrm{I}_{2,1}
	:=
	\int_{\Omega_{T,R}}
	\frac{
		\|B^{-1}(\sigma(s,x)-\sigma(s,y))\|^2
	}{
		|B^{-1}(x-y)|^2
	}
	\,\pi_s(\d x,\d y)\,\d s.
	$$
	Applying Lemma \ref{A variant maximal function estimate} (see also Remark \ref{rem-assump1}(iii)), we have for a.e. $t\in[0,T]$ and for a.e. $x,y\in\Rd$,
	\begin{align*}
		&\quad\frac{
			\|B^{-1}(\sigma(t,x)-\sigma(t,y))\|^2
		}{
			|B^{-1}(x-y)|^2
		}\\
		&\lesssim
		\frac{
			\|B^{-1}(\sigma_1(t,x)-\sigma_1(t,y))\|^2
		}{
			|B^{-1}(x-y)|^2
		}
		+
		\frac{
			\|B^{-1}(\sigma_2(t,x)-\sigma_2(t,y))\|^2
		}{
			|B^{-1}(x-y)|^2
		}
		\\
		&\lesssim
		\frac{1}{\delta_1^2}
		\left(\widetilde{\M}(\nabla_{x_1}\sigma_1)\right)^2(t,x)
		+
		\frac{1}{\delta_1^2}
		\left(\widetilde{\M}(\nabla_{x_1}\sigma_1)\right)^2(t,y)
		\\
		&\quad
		+
		\frac{1}{\delta_2^2}
		\left(\widetilde{\M}(\nabla\sigma_2)\right)^2(t,x)
		+
		\frac{1}{\delta_2^2}
		\left(\widetilde{\M}(\nabla\sigma_2)\right)^2(t,y).
	\end{align*}
	Combining the fact that $\rho,\hat{\rho}\in L^\infty_t(L^1_x\cap L^\infty_x)$, assumption $(\mathcal{H}^{\sigma})$, Lemma \ref{Operator bounds}, by the similar argument as above, we obtain
	$$\mathrm{I}_{2,1}\lesssim 1.$$
	Taking $\alpha=1$,
	\begin{align}\label{pf-thm partial(2)-4}
		\mathrm{I}_{2}\lesssim \frac{1}{\delta_1^2}\|\sigma-\hat{\sigma}\|^2_{L^2([0,T];L^{2r}(B_R))}+1.
	\end{align}
	With $\alpha=1$, by \eqref{pf-thm partial(2)-1}, we have
	\begin{align}\label{pf-thm partial(2)-5}
		\E\log\left(1+\sqrt{|B^{-1}(X_0-\hat X_0)|^2+\alpha}\right)
		&\leq
		\E\log\left(1+|B^{-1}(X_0-\hat X_0)|\right)+ 1\cr
		&=\D_{\delta_1,\delta_2}(\nu,\hat\nu)+ 1.
	\end{align}
	
	Finally, combining the estimates  \eqref{pf-thm partial(2)-2},  \eqref{pf-thm partial(2)-3}, \eqref{pf-thm partial(2)-4} and  \eqref{pf-thm partial(2)-5} and taking the supremum on $[0,T]$, we complete the proof.
\end{proof}

\subsection{Rates of convergence}\label{rate}\hskip\parindent
The goal of this subsection is to show that the distance $\D_{\delta_1,\delta_2}$ plays a crucial role in deriving explicit weak convergence rates for solutions to
the FPE, which further illustrates its advantage over the classical $L^1$-Wasserstein distance. To make this point precise, we specialize Assumption \ref{assump1} to the case $b_1=0$ and $\sigma=0$, with $b_2\in L^p([0,T]\times\Rd; \R^{d_2})$, so that the integrability now holds over the whole
space rather than only locally.

\textbf{Perturbations of $b_2$.}  We consider a sequence $b_2^n\to b_2$ in $L_t^1L_x^1$ as $n\to\infty$, where each $b_2^n$ satisfies the conditions in Assumption~\ref{assump1},
while keeping
$$\rho_0^n=\rho_0,\qquad b_1^n=b_1=0, \qquad \sigma^n=\sigma=0.$$

For each $n\geq1$, let
$\rho^n,\rho\in L_t^\infty(L_x^1\cap L_x^\infty)$
denote the corresponding weak solutions to the FPE with drift components
$b_2^n$ and $b_2$, respectively. Set
$$
e_n:=\|b_2^n-b_2\|_{L_t^1L_x^1},
$$
so that $e_n\to0$ as $n\to\infty$. By Theorem~\ref{thm partial}(2), under the
present specialization, the right-hand side of the stability estimate is
independent of $R$. Hence, letting $R\to\infty$, we obtain
$$
\sup_{t\in[0,T]}
\D_{\delta_1,\delta_2}
\bigl(\rho^n(t,\cdot),\rho(t,\cdot)\bigr)
\lesssim
\frac{e_n}{\delta_1}+1
$$
for every admissible pair $(\delta_1,\delta_2)$.

To quantify the weak convergence, we choose vanishing
admissible
parameters $(\delta_1^n,\delta_2^n)$ such that the corresponding Kantorovich--Rubinstein distance remains uniformly bounded. By the parameter choice in the proof of Theorem~\ref{thm partial}(2), for all sufficiently large $n$, we may take
$$
\delta_1^n:=e_n,
\quad
\delta_2^n:=(\delta_1^n)^{\frac{s}{s+1}},
\qquad s\in(0,1),
$$
and therefore,
$$
\sup_{t\in[0,T]}
\D_{\delta_1^n,\delta_2^n}
\bigl(\rho^n(t,\cdot),\rho(t,\cdot)\bigr)
\lesssim1,
$$
with a constant independent of $n$.

We now show quantitatively that the shrinking of these scales implies weak
convergence. Fix $0<\theta<1$ and set
$$
\eta_i^n:=(\delta_i^n)^\theta,\qquad i=1,2.
$$
Then for each $i=1,2$, $\eta_i^n\to0$ and $\eta_i^n/\delta_i^n\to\infty$ as $n\rightarrow\infty$. For each
$t\in[0,T]$, let $\pi_{n,t}$ be an optimal coupling associated with
$\D_{\delta_1^n,\delta_2^n}
(\rho^n(t,\cdot),\rho(t,\cdot))$, and let $B_n$ be the matrix in
\eqref{diag-B} corresponding to $(\delta_1^n,\delta_2^n)$. Writing
$x=(x_1,x_2)\in\R^{d_1}\times\R^{d_2}$ and $y=(y_1,y_2)\in\R^{d_1}\times\R^{d_2}$, define
$$
A_{i,n}:=\{(x,y):\ |x_i-y_i|>\eta_i^n\},
\qquad i=1,2.
$$
Then
$$
\pi_{n,t}(A_{i,n})
\log\left(1+\frac{\eta_i^n}{\delta_i^n}\right)
\leq
\D_{\delta_1^n,\delta_2^n}
\bigl(\rho^n(t,\cdot),\rho(t,\cdot)\bigr)
\lesssim1,\qquad i=1,2,
$$
and hence
$$
\sup_{t\in[0,T]}\pi_{n,t}(A_{1,n}\cup A_{2,n})
\lesssim
\frac{1}{\log(1+(\delta_1^n)^{\theta-1})}
+
\frac{1}{\log(1+(\delta_2^n)^{\theta-1})}.
$$
Consequently, for every bounded Lipschitz function $f$ on $\Rd$, using the coupling $\pi_{n,t}$ and splitting the integral over $(A_{1,n}\cup A_{2,n})^c$ and $A_{1,n}\cup A_{2,n}$, we obtain
$$
\begin{aligned}
	&\sup_{t\in[0,T]}
	\left|
	\int_{\Rd} f(x)\rho^n(t,x)\,\d x
	-
	\int_{\Rd} f(x)\rho(t,x)\,\d x
	\right|
	\\
	&\leq
	\sup_{t\in[0,T]}
	\int_{\Rd\times\Rd}
	|f(x)-f(y)|\,\pi_{n,t}(dx,dy)
	\\
	&\leq
	\operatorname{Lip}(f)(\eta_1^n+\eta_2^n)
	+
	2\|f\|_{L^\infty(\Rd)}
	\sup_{t\in[0,T]}
	\pi_{n,t}(A_{1,n}\cup A_{2,n})
	\\
	&\lesssim
	\operatorname{Lip}(f)\left(e_n^\theta + e_n^{\frac{s\theta}{s+1}}\right)
	+
	2\|f\|_{L^\infty(\Rd)}\bigg[\frac{1}{\log\Big(1+e_n^{-(1-\theta)}\Big)}
	+
	\frac{1}{
		\log\Big(
		1+e_n^{-\frac{s(1-\theta)}{s+1}}
		\Big)}\bigg]\\
	&\longrightarrow0,\qquad\mbox{as }n\rightarrow\infty,
\end{aligned}
$$
where $\operatorname{Lip}(f)$ stands for the Lipschitz constant of $f$. Thus, $\rho^n$ converges weakly to $\rho$, uniformly in $t\in[0,T]$, and the preceding estimate gives an explicit quantitative rate of convergence against bounded Lipschitz test functions.

\textbf{Perturbations of the initial data.}
We now assume $\rho_0^n\to\rho_0$ in $L_x^1\cap L_x^\infty$ as $n\rightarrow\infty$, while
$b^n=b,\ \sigma^n=\sigma$ for $n\geq1$.
Combining the linearity of the FPE with the energy estimates, we obtain, for every
$1\leq p\leq\infty$,
$$
\|\rho^n-\rho\|_{L_t^\infty L_x^p}
\lesssim
\|\rho_0^n-\rho_0\|_{L_x^p}.
$$
Consequently, $\rho^n$ converges strongly to $\rho$ in $L_t^\infty L_x^p$ for every
$1\leq p\leq\infty$, with convergence rate $\|\rho_0^n-\rho_0\|_{L_x^p}$. This rate
is clearly optimal.

\section{Gradient of drift fields as singular integrals}\label{singular}
\subsection{Main result}\hskip\parindent
We introduce several classes of singular kernels of convolution type, and recall the basic properties
of the associated singular integral operators; see \cite{Grafakos2014a} for more details.
For a function $K$, we write $\widehat{K}$ for its Fourier transform whenever it is well defined.

\begin{definition}[Singular kernel]\label{def:singular-kernel}
	Let $K\in\mathscr{S}'(\R^d)$ which coincides with a locally integrable function on $\R^d\setminus\{0\}$. We say that
	$K$ is a singular kernel on $\Rd$ if there exist constants $A_1,A_2\geq0$ such that
	\begin{align*}
		&|\widehat{K}(\xi)|\leq A_1,\quad \xi\in\R^d,\\
		&\int_{|x|\geq 2|y|}|K(x-y)-K(x)|\,\d x\leq A_2,\quad 
		y\in \R^d\setminus\{0\}.
	\end{align*}
\end{definition}

The Calder\'{o}n--Zygmund theory yields the following results.
\begin{theorem}[Calder\'{o}n--Zygmund]\label{thm:calderon-zygmund}
	Let $K$ be a singular kernel as in Definition \ref{def:singular-kernel}. For $u\in L^2(\Rd)$,  
	define the operator associated with $K$ as
	$$
	Su=K*u,
	$$
	in the sense of multiplication in the Fourier variable. Then, for every
	$1<p<\infty$, there exists a constant $C_{d,p}>0$ such that
	\begin{equation}\label{Lp estimate}
		\|Su\|_{L^p(\Rd)}
		\leq C_{d,p}
		(A_1+A_2)
		\|u\|_{L^p(\Rd)},
		\qquad u\in L^p(\Rd)\cap L^2(\Rd).
	\end{equation}
	For $p=1$, there exists a constant $C_d>0$ such that one has the weak-type estimate
	\begin{equation}\label{L1 estimate}
		\|Su\|_{L^{1,\infty}(\Rd)}
		\leq C_d
		(A_1+A_2)
		\|u\|_{L^1(\Rd)},
		\qquad u\in L^1(\Rd)\cap L^2(\Rd).
	\end{equation}
	Consequently, $S$ extends uniquely to a bounded operator on $L^p(\R^d)$ for $1<p<\infty$, and to a bounded operator from $L^1(\R^d)$ to $L^{1,\infty}(\R^d)$, still denoted by $S$, for which the above estimates hold for all $u$ in the respective spaces.
\end{theorem}
The operator obtained in this way is called the singular integral operator associated with $K$. For $p=1$, the extension
$S:  L^1(\Rd)\to L^{1,\infty}(\Rd)$
may be
misleading, since a function in $L^{1,\infty}(\Rd)$ is in general not locally
integrable and therefore need not define a distribution.  To circumvent this, we introduce a different extension into the tempered distributions $S^{\rm dis}: L^1(\Rd)\to \mathscr{S}'(\Rd)$,  defined via the duality
\begin{equation}\label{distribution extension}
	\<S^{\rm dis} u,\varphi\>:=\<u,\widetilde{S}\varphi\>,
	\qquad \varphi\in \mathscr{S}(\Rd),
\end{equation}
where $\widetilde{S}$ is the singular integral operator associated with
$\widetilde{K}(x):=K(-x)$. It is important to emphasize that this distributional extension is not the same as the extensions $S:  L^1(\Rd)\to L^{1,\infty}(\Rd)$. The two  should not be identified.

For our purposes, we use a more regular class of kernels, namely, singular kernels of
fundamental type.
\begin{definition}[Singular kernel of fundamental type]\label{def:singular-kernel-fundamental-type}
	We say that $K: \R^d\setminus\{0\}\rightarrow\R$ is a singular kernel of fundamental type if it coincides with a $C^1$ function on $\mathbb{R}^d\setminus\{0\}$ and there exist constants $C_1,C_2\geq 0$ such that
	\begin{itemize}
		\item[(i)](Size and gradient condition)
		$$| K(x)|+ |x||\nabla K(x)|\leq \frac{C_1}{|x|^{d}},\qquad  x\neq 0;$$
		
		\item[(ii)](Cancellation condition)
		\[
		\sup_{0<R_1<R_2<\infty}
		\left|
		\int_{B_{R_2}\setminus B_{R_1}}K(x)\,\mathrm{d}x
		\right|
		\leq C_2.
		\]
	\end{itemize}
	As a consequence, $K$ admits an extension to a tempered distribution on $\R^d$, which is a singular kernel and is unique up to adding a constant multiple of the Dirac mass at the origin.
\end{definition}

We also recall
the following standard notion
from the distribution theory.
\begin{definition}
	For $f\in \mathscr{D}'(\Rd)$, its distributional derivative $\partial_j f$ is
	defined by
	$$
	\<\partial_j f,\psi\>:=-\<f,\partial_j\psi\>,
	\qquad \psi\in \mathscr{D}(\Rd).
	$$
\end{definition}

For the statements of our main results, we require the following assumption.
\begin{assumption}\label{assump stability}
	Let $K=(K_{ij})_{1\leq i,j\leq d}$ be a matrix-valued kernel such that each $K_{ij}$ is a singular kernel of fundamental type, and let $w=(w_i)_{1\leq i\leq d}$ such that each  $w_i\in L^1([0,T];L^1(\R^d))$.
	\begin{itemize}
		\item[]$(\widetilde{\mathcal{H}}^{\nabla b})$  $\nabla b = K*w=(K_{ij}*w_i)_{1\leq i,j\leq d}$, and the convolutions are understood in the distributional sense of \eqref{distribution extension}.
		
		\item[]$(\widetilde{\mathcal{H}}^{b})$
		$b \in L^{p,\infty}([0,T]\times\mathbb{R}^d;\R^d)$ for some $p>1$.
		
		\item[]$(\widetilde{\mathcal{H}}^{\sigma})$
		$\sigma \in L^2([0,T]; L^\infty(\mathbb{R}^d; \mathbb{R}^{d\times d}))\text{ with }\left( \sum_{i,j,k=1}^d |\partial_k \sigma_{ij}|^2 \right)^{1/2} \in L^2([0,T]; L^{2r}(\mathbb{R}^d))$ for some $r>1$.
	\end{itemize}
\end{assumption}
We note that assumption $(\widetilde{\mathcal{H}}^{b})$ implies $b\in L^1_{\rm loc}([0,T]\times \R^d;\R^d)$.  This follows immediately by applying the interpolation inequality \eqref{interp 1} with $r=1$ and $\ell=p>1$ on any compact subset of $[0,T]\times \R^d$. Moreover, under assumption $(\widetilde{\mathcal{H}}^{\sigma})$, we have $\sigma(t,\cdot)\in W^{1,1}_{\rm loc}(\R^d; \R^{d\times d})$ for almost every $t\in[0,T]$.

As before, we consider an increasing logarithmical
cost, which generates the Kantorovich--Rubinstein distance
$$\D_{\delta}(\mu, \nu):=\inf_{\pi \in \mathscr{C}(\mu, \nu)}\int_{\Rd\times \Rd} \log\left(1+\frac{|x-y|}{\delta}\right) \,\pi(\d x, \d y), \qquad \delta>0,\  \mu, \nu \in \P(\Rd).$$

We now state the stability estimate.
\begin{theorem}\label{stability estimate 2}
	Let $(b,\sigma,\rho_0)$ and $(\hat b,\hat\sigma,\hat\rho_0)$ be two triples
	satisfying Assumption \ref{assump stability} and
	$(\mathcal{H}^{\nu})$. Let
	$\rho,\hat\rho\in L^\infty([0,T]; L^1(\R^d)\cap L^\infty(\R^d))$ be the corresponding weak
	solutions to the FPE \eqref{FPE}. Then, for every $\e>0$, there exists a
	constant $C_\e>0$ such that, for every $\delta>0$,
	\begin{align*}
		\sup_{t\in [0,T]}	\D_\delta(\rho(t,\cdot),\hat\rho(t,\cdot))
		&\lesssim
		\D_\delta(\rho_0,\hat\rho_0)
		+\frac{1}{\delta}
		\int_0^T\int_{\Rd}
		|b(t,x)-\hat b(t,x)|\rho(t,x)\,\d x\d t  \\
		&\quad
		+\frac{1}{\delta^2}
		\int_0^T\int_{\Rd}
		\|\sigma(t,x)-\hat\sigma(t,x)\|^2\rho(t,x)\,\d x\d t \\
		&\quad +\e\left[1+\log\left(\frac{1}{\e\delta}\right)\right]+C_\e.
	\end{align*}
	Here the implicit constant in $\lesssim$ is independent of $\e$ and $\delta$.
\end{theorem}

\begin{remark}
	In Theorem \ref{stability estimate 2}, the assumption $(\widetilde{\mathcal{H}}^{\sigma})$ is imposed to include the identity matrix in our framework. It can be replaced by the standard Sobolev regularity condition
	$$ \sigma\in L^2([0,T];W^{1,2r}(\R^{d};\R^{d\times d}))$$
	for some $r>1$.
\end{remark}

The stability estimate in Theorem \ref{stability estimate 2} can be restated as follows.
\begin{corollary}\label{cor-stability2}
	Under the assumptions of Theorem \ref{stability estimate 2}, for every $\e>0$, there exists a
	constant $C_\e>0$ such that, for every $\delta>0$,
	\begin{align*}
		\sup_{t\in [0,T]}	\D_\delta(\rho(t,\cdot),\hat\rho(t,\cdot))
		&\lesssim
		\D_\delta(\rho_0,\hat\rho_0)
		+\frac{1}{\delta^2}\|\sigma-\hat{\sigma}\|_{L^2([0,T]; L^\infty(\R^d))}^2 \|\rho\|_{L^\infty([0,T];L^1(\R^d))}\\
		&\quad +\frac{1}{\delta}\|b-\hat{b}\|_{L^{p,\infty}([0,T]\times\R^d)}  \|\rho\|_{L^1([0,T]\times\R^d)}^{1/p'}\|\rho\|_{L^\infty([0,T]\times\R^d)}^{1/p}\\
		&\quad +\e\left[1+\log\left(\frac{1}{\e\delta}\right)\right]+C_\e,
	\end{align*}
	where $p'$ is the conjugate exponent of $p$, and the implicit constant in $\lesssim$ is independent of $\e$ and $\delta$.
\end{corollary}
\begin{proof}
	Let $E=[0,T]\times\R^d$. Since $\rho\in L^\infty([0,T]; L^1(\R^d)\cap L^\infty(\R^d))$, choosing a product-measurable representative and applying Tonelli's theorem immediately yield $\rho\in L^1(E)\cap L^\infty(E)$. By standard real interpolation, it holds that
	$$\left(L^1(E),L^\infty(E)\right)_{1/q,1}=L^{q',1}(E),\qquad q\in(1,\infty),\,q':=\frac{q}{q-1}$$
	with an equivalent norm. Hence, $\rho\in L^{q,1}(E)$ for every $q\in(1,\infty)$, and there exists a constant $C_q>0$ such that
	\begin{align*}\label{pf-cor-stability2-1}
		\|\rho\|_{L^{q,1}(E)}\leq C_q\|\rho\|_{L^1(E)}^{1/q}\|\rho\|_{L^\infty(E)}^{1/q'}.
	\end{align*}
	Combining this inequality with \cite[Exercise 1.4.1(b)]{Grafakos2014a}, we obtain for any $p\in(1,\infty)$ with $p'=\frac{p}{p-1}$,
	\begin{align*}
		\int_0^T\int_{\R^d}|b(t,x)-\hat b(t,x)|\rho(t,x)\,\d x\d t&\leq \|b-\hat{b}\|_{L^{p,\infty}(E)}\|\rho\|_{L^{p',1}(E)}\cr
		&\lesssim \|b-\hat{b}\|_{L^{p,\infty}(E)}\|\rho\|_{L^1(E)}^{1/p'}\|\rho\|_{L^\infty(E)}^{1/p}.
	\end{align*}
	By the assumption $\sigma, \hat{\sigma} \in L^2([0,T]; L^\infty(\mathbb{R}^d; \mathbb{R}^{d\times d}))$, a direct estimate gives
	\begin{align*}
		\int_0^T\int_{\Rd}\|\sigma(t,x)-\hat\sigma(t,x)\|^2\rho(t,x)\,\d x\d t\leq \|\sigma-\hat{\sigma}\|_{L^2([0,T]; L^\infty(\R^d))}^2 \|\rho\|_{L^\infty([0,T];L^1(\R^d))}.
	\end{align*}
\end{proof}

As an immediate consequence of Theorem \ref{stability estimate 2}, we have the following uniqueness.
\begin{corollary}\label{cor:uniqueness}
	Under the assumptions of Theorem \ref{stability estimate 2}, weak solutions to
	the FPE \eqref{FPE} in the class $L^\infty([0,T]; L^1(\R^d)\cap L^\infty(\R^d))$
	are unique.
\end{corollary}
\begin{proof}
	Let $b=\hat{b}$, $\sigma=\hat{\sigma}$ and $\rho_0=\hat\rho_0$. By Theorem \ref{stability estimate 2},
	$$\D_{\delta}(\rho(t,\cdot), \hat{\rho}(t,\cdot)) \lesssim \e\left[1+\log\left(\frac{1}{\e\delta}\right)\right]+C_{\e},\quad t\in[0,T],\,\epsilon>0.
	$$
	If $0<\delta<e^{-1}$, that is, $|\log \delta|>1$, then
	$$\frac{1}{|\log \delta|}\left[1+\log\left(\frac{1}{\e\delta}\right)\right] \leq\frac{1+|\log \delta|+|\log \e|}{|\log \delta|}\leq 2+|\log \e|.$$
	Hence, for an arbitrarily small constant $a>0$, there exists $0<\e_0=\e_0(a)<1$ such that
	$$\frac{\e_0}{|\log \delta|}\left[1+\log\left(\frac{1}{\e_0\delta}\right)\right] \leq  \frac{a}{2}.$$
	For this fixed $\e_0>0$, we may choose $\delta\in(0,e^{-1})$ small enough that
	$$\frac{
		C_{\e_0}
	}{|\log \delta|}\leq \frac{a}{2}.$$
	Thus, for every $t\in[0,T]$,
	\begin{equation}\label{log lim}
		\lim_{\delta\to 0} \frac{\D_{\delta}(\rho(t,\cdot), \hat{\rho}(t,\cdot))}{|\log \delta|}=0.
	\end{equation}
	Remark 2.1(2)
	in \cite{Li 2019} shows that \eqref{log lim} implies $\rho=\hat{\rho}$ in $L_t^\infty(L_x^1\cap L_x^\infty)$.
\end{proof}

\subsection{Proof of Theorem \ref{stability estimate 2}}\hskip\parindent
Compared with the previous section, the primary challenge in proving Theorem \ref{stability estimate 2} is that the Morrey-type
estimate from Lemma \ref{Morrey-type estimate} no longer applies to the relevant
difference quotients for the drift field $b$ introduced in Assumption \ref{assump stability}. We therefore follow \cite{Bouchut 2013} and introduce a
smooth variant of the usual maximal function. Given $\varphi\in C_c^\infty(\Rd)$, for $\epsilon>0$, set
$$\varphi^{\e}(x):=\frac{1}{\e^d}\varphi\left(\frac{x}{\e}\right), \quad x\in \Rd.$$ For each component, define
$$\M_{\varphi}(K_{ij}*w_i)(t,x):=\sup_{\e>0}\left|\<(K_{ij}*w_i)(t,\cdot),\varphi^\e(x-\cdot)\>\right|,\quad (t,x)\in [0,T]\times \Rd,$$
in the sense of \eqref{distribution extension}, and set
$$
\M_\varphi(K*w)(t,x)
:=
\left(
\sum_{i,j=1}^d
\bigl|\M_\varphi(K_{ij}*w_i)(t,x)\bigr|^2
\right)^{1/2},
\quad (t,x)\in [0,T]\times \Rd.
$$

The first part of the next result is from \cite[Proposition~4.2]{Bouchut 2013}, and the second part is an easy consequence of \cite[Theorem 3.3]{Bouchut 2013} and the Marcinkiewicz interpolation theorem (cf. \cite[Theorem 1.3.2]{Grafakos2014a}) since the operator $\M_\varphi$ is sublinear.
\begin{lemma}\label{lem:weighted-maximal}
	Let $b\in L_{\mathrm{loc}}^1([0,T]\times\Rd;\R^d)$ and assume that
	$(\widetilde{\mathcal{H}}^{\nabla b})$
	holds. Then there exist a constant $C_d>0$ and
	a function $\varphi\in C_c^\infty(\Rd)$ such that, for a.e. $t\in[0,T]$, there
	is an $\mathcal{L}^d$-negligible set $N_t\subset\Rd$ with
	\begin{equation}\label{crippa morrey}
		|b(t,x)-b(t,y)|
		\leq
		C_d |x-y|
		\left(
		\M_\varphi(\nabla b)(t,x)
		+
		\M_\varphi(\nabla b)(t,y)
		\right)
	\end{equation}
	for all $x,y\notin N_t$. 
	Moreover, for every $1<p\leq2$, there exists a constant $C_{d,p}>0$ such that, for a.e. $t\in[0,T]$,
	\begin{align}\label{crippa maximal estimates}
		\|\M_\varphi(\nabla b)(t,\cdot)\|_{L^p(\Rd)}\leq C_{d,p}\|w(t,\cdot)\|_{L^p(\R^d)},
	\end{align}
	and for $p=1$, there exists a constant $C'_{d}>0$ such that, for a.e. $t\in[0,T]$,
	\begin{align}\label{crippa maximal estimates+}
		\|\M_\varphi(\nabla b)(t,\cdot)\|_{L^{1,\infty}(\Rd)}\leq C'_d\|w(t,\cdot)\|_{L^1(\R^d)}.
	\end{align}
\end{lemma}

For the application below, we use the following variant of \eqref{crippa morrey}.
\begin{lemma}\label{lem:modified-crippa-morrey}
	Under the assumption of Lemma \ref{lem:weighted-maximal}, for any $\epsilon>0$, there exist measurable functions
	$G_\e^1,G_\e^2:[0,T]\times\Rd\to \R_+$ and a constant $C_\e>0$ such that
	$$
	\|G_\e^1\|_{L^1([0,T];L^{1,\infty}(\Rd))}\leq \e,
	\qquad
	\|G_\e^2\|_{L^1([0,T];L^2(\Rd))}\leq C_\e,
	$$
	and moreover, for a.e. $t\in[0,T]$, there exists an $\mathcal{L}^d$-negligible set
	$N_t\subset\Rd$ such that for all $x,y\notin N_t$,
	\begin{equation}\label{b gradient morrey}
		|b(t,x)-b(t,y)|
		\leq C_d
		|x-y|
		\left(	G_\e^1(t,x)+G_\e^1(t,y)+G_\e^2(t,x)+G_\e^2(t,y) \right),
	\end{equation}
	where $C_d>0$ is a constant depending only on $d$.
\end{lemma}

\begin{proof}
	Fix $\e>0$. Since $w\in L^1([0,T];L^1(\R^d;\R^d))$, by
	a standard truncation, write $w=w_1^\eta+w_2^\eta$ with
	$$
	\|w_1^\eta\|_{L_t^1L_x^1}\leq \eta,
	\qquad
	\|w_2^\eta\|_{L_t^1L_x^2}\leq C_\eta ,
	$$
	where $\eta>0$ depends on $\epsilon$ and is chosen below. For instance, we may choose $w_1^\eta=w\mathbbm{1}_{\{|w|>R\}}$ for large enough $R>0$ to ensure that $\|w_1^\eta\|_{L_t^1L_x^1}\leq \eta$, and hence, $w_2^\eta=w\mathbbm{1}_{\{|w|\leq R\}}$ clearly satisfies the desired estimate. Define
	\begin{align}\label{M_varphi}
		G_\e^1:=\M_\varphi(K*w_1^\eta),
		\qquad
		G_\e^2:=\M_\varphi(K*w_2^\eta),
	\end{align}
	where  $\varphi\in C_c^\infty(\Rd)$.
	Then \eqref{b gradient morrey} follows from \eqref{crippa morrey} and the
	sublinearity of $\M_\varphi$. Finally, by the weak $(1,1)$ estimate \eqref{crippa maximal estimates+} and the
	$L^2$ boundedness \eqref{crippa maximal estimates}, we may choose $\eta>0$
	sufficiently small so that the stated bounds on $G_\e^1$ and $G_\e^2$ hold.
\end{proof}

\begin{proof}[Proof of Theorem \ref{stability estimate 2}]
	As in the proof of Theorem \ref{thm partial}(2), we reduce to the SDEs via the superposition principle in Theorem \ref{sp}, where the required integrability is guaranteed by our assumptions. Lemma \ref{lem-common} then allows us to assume, without loss of generality, that the two SDE solutions are coupled on a common probability space with the same Brownian motion, and with an optimal initial coupling. So, we consider the coupled SDEs driven by the same Brownian motion $(W_t)_{t\ge0}$:
	$$
	\d X_t=b(t,X_t)\d t+\sigma(t,X_t)\d W_t,
	\qquad \L_{X_0}=\nu,
	$$
	and
	$$
	\d \hat{X}_t=\hat{b}(t,\hat{X}_t)\d t+\hat{\sigma}(t,\hat{X}_t)\d W_t,
	\qquad \L_{\hat{X}_0}=\hat{\nu},
	$$
	where $\nu=\rho_0\mathcal{L}^d$ and $\hat{\nu}=\hat{\rho}_0\mathcal{L}^d$. The initial coupling is chosen optimally with respect to the Kantorovich--Rubinstein distance, i.e.,
	\begin{align}\label{pf-se2-1}
		\E \log \left(1+\frac{\sqrt{|X_0-\hat{X}_0|^2}}{\delta}\right)=\mathcal{D}_{\delta}(\nu,\hat{\nu}).
	\end{align}By the superposition principle again, the time marginals satisfy
	\[
	\mathcal{L}_{X_t}(dx)=\rho(t,x)dx,\quad \mathcal{L}_{\hat X_t}(dx)=\hat\rho(t,x)dx,\qquad t\in[0,T].
	\]

	Let $\alpha=\delta^2>0$. By It\^o's formula,
	\begin{align}\label{pf-1}
		&\E \log \left(1+\frac{\sqrt{|X_t-\hat{X}_t|^2+\alpha}}{\delta}\right)\cr
		&\leq
		\E \log \left(1+\frac{\sqrt{|X_0-\hat{X}_0|^2+\alpha}}{\delta}\right)
		+\int_0^T\int_{\Rd\times\Rd}
		\frac{|b(t,x)-\hat{b}(t,y)|}{\delta+|x-y|}
		\,\pi_t(\d x,\d y)\d t \cr
		&\quad
		+\int_0^T\int_{\Rd\times\Rd}
		\frac{\|\sigma(t,x)-\hat{\sigma}(t,y)\|^2}
		{\delta\sqrt{|x-y|^2+\alpha}+|x-y|^2+\alpha}
		\,\pi_t(\d x,\d y)\d t,\quad t\in[0,T],
	\end{align}
	where $\pi_t=\L_{(X_t,\hat{X}_t)}$ for $t\in[0,T]$. Set
	$$
	\mathrm{J}_1
	:=
	\int_0^T\int_{\Rd\times\Rd}
	\frac{|b(t,x)-\hat{b}(t,y)|}{\delta+|x-y|}
	\,\pi_t(\d x,\d y)\d t,
	$$
	and
	$$
	\mathrm{J}_2
	:=
	\int_0^T\int_{\Rd\times\Rd}
	\frac{\|\sigma(t,x)-\hat{\sigma}(t,y)\|^2}
	{\delta\sqrt{|x-y|^2+\alpha}+|x-y|^2+\alpha}
	\,\pi_t(\d x,\d y)\d t.
	$$
	
	We first estimate $\mathrm{J}_2$. By the Morrey-type estimate in
	Lemma \ref{Morrey-type estimate},  the triangle inequality and H\"{o}lder's inequality,
	\begin{align}\label{pf-J2}
		\mathrm{J}_2
		&\leq
		\int_0^T\int_{\Rd\times\Rd}
		\left\{
		\frac{\|\sigma(t,x)-\hat{\sigma}(t,x)\|^2}{\delta\sqrt{\alpha}}
		+
		\frac{\|\hat{\sigma}(t,x)-\hat{\sigma}(t,y)\|^2}{|x-y|^2}
		\right\}
		\,\pi_t(\d x,\d y)\d t \cr
		&\lesssim
		\frac{1}{\delta^2}
		\int_0^T\int_{\Rd}
		\|\sigma(t,x)-\hat{\sigma}(t,x)\|^2\rho(t,x)\,\d x\d t \cr
		&\quad
		+
		\int_0^T
		\left(
		\int_{\Rd} |\M(\nabla\hat{\sigma})(t,x)|^{2r}\,\d x
		\right)^{1/r}
		\left(
		\|\rho(t,\cdot)\|_{L_x^{r'}}
		+
		\|\hat{\rho}(t,\cdot)\|_{L_x^{r'}}
		\right)\d t \cr
		&\lesssim
		\frac{1}{\delta^2}
		\int_0^T\int_{\Rd}
		\|\sigma(t,x)-\hat{\sigma}(t,x)\|^2\rho(t,x)\,\d x\d t
		+1,
	\end{align}
	where $r'$ is the conjugate number of $r$, and the last inequality is due to the assumption that $\rho,\hat{\rho}\in L^\infty_t(L^1_x\cap L^\infty_x)$ and
	$( \sum_{i,j,k=1}^d |\partial_k \hat{\sigma}_{ij}|^2 )^{1/2} \in L^2_tL^{2r}_x$, and the strong $(2r,2r)$ boundedness of $\mathcal{M}$.
	
	We next estimate $\mathrm{J}_1$. By the triangle inequality,
	\begin{align}\label{pf-J1-1}
		\mathrm{J}_1&\leq\frac{1}{\delta}\int_0^T\int_{\Rd}|b(t,x)-\hat{b}(t,x)|\rho(t,x)\,\d x\d t  \cr
		&\quad +
		\int_0^T\int_{\Rd\times\Rd}\frac{|\hat{b}(t,x)-\hat{b}(t,y)|}{\delta+|x-y|}\,\pi_t(\d x,\d y)\d t \cr
		&=:\frac{1}{\delta}\int_0^T\int_{\Rd}|b(t,x)-\hat{b}(t,x)|\rho(t,x)\,\d x\d t   +   \mathrm{J}_{1,1}.
	\end{align}
	By \eqref{b gradient morrey} with $b$ replaced by $\hat{b}$, for any $\epsilon>0$,
	\begin{align*}
		&\quad\frac{|\hat{b}(t,x)-\hat{b}(t,y)|}{\delta+|x-y|}\\
		&\lesssim
		\min\left\{
		\frac{|\hat{b}(t,x)-\hat{b}(t,y)|}{\delta},
		G_\e^1(t,x)+G_\e^1(t,y)
		\right\}
		+G_\e^2(t,x)+G_\e^2(t,y) \\
		&\lesssim
		\min\left\{\frac{|\hat{b}(t,x)|}{\delta},G_\e^1(t,x)\right\}
		+
		\min\left\{\frac{|\hat{b}(t,y)|}{\delta},G_\e^1(t,y)\right\} \\
		&\quad
		+
		\min\left\{\frac{|\hat{b}(t,y)|}{\delta},G_\e^1(t,x)\right\}
		+
		\min\left\{\frac{|\hat{b}(t,x)|}{\delta},G_\e^1(t,y)\right\} \\
		&\quad
		+G_\e^2(t,x)+G_\e^2(t,y),
	\end{align*}
	where the functions $G_\e^1$ and $G_\e^2$ are defined similarly as \eqref{M_varphi}.
	It is therefore enough to estimate the representative terms.
	
	Set
	$$
	\eta_1(t,x)
	:=
	\min\left\{\frac{|\hat{b}(t,x)|}{\delta},G_\e^1(t,x)\right\},
	\quad
	\eta_2(t,x,y)
	:=
	\min\left\{\frac{|\hat{b}(t,y)|}{\delta},G_\e^1(t,x)\right\},
	$$
	and
	$$
	\eta_3(t,x):=G_\e^2(t,x).
	$$
	Define a finite measure $\mu$ on $\Omega_T:=[0,T]\times\Rd\times\Rd$ by
	$$
	\mu(\d t,\d x,\d y)
	=
	\mathbbm{1}_{[0,T]}(t)\pi_t(\d x,\d y)\d t .
	$$

	Since $\rho,\hat{\rho}\in L_t^\infty(L_x^1\cap L_x^\infty)$, by Lemma \ref{lem:modified-crippa-morrey},
	\begin{align}\label{pf-eta3}
		\|\eta_3\|_{L^1(\Omega_T,\mu)}
		\lesssim
		\|G_\e^2\|_{L_t^1L_x^2}
		\leq C_\e .
	\end{align}
	
	We now estimate $\eta_1$ and $\eta_2$. Since $\eta_1\leq G_\e^1$ and $\rho\in L_t^\infty(L_x^1\cap L_x^\infty)$, by Lemma \ref{lem:modified-crippa-morrey},
	\begin{equation*}
		\begin{split}
			\|\eta_1\|_{L^{1,\infty}(\Omega_T,\mu)}
			&=
			\sup_{\lambda>0}
			\lambda\,
			\mu\left(\{(t,x,y)\in\Omega_T:\ 
			|\eta_1(t,x)|
			>\lambda\}\right) \\
			&\leq
			\sup_{\lambda>0}
			\lambda
			\int_0^T
			\int_{\{x:\ |G_\e^1(t,x)|>\lambda\}}
			\rho(t,x)\,\d x\d t \\
			&\leq
			\|\rho\|_{L^\infty([0,T]\times\Rd)}
			\sup_{\lambda>0}
			\lambda
			\int_0^T
			\mathcal{L}^d(\{x:\ |G_\e^1(t,x)|>\lambda\})\,\d t \\
			&\lesssim
			\|G_\e^1\|_{
				L^1([0,T]; L^{1,\infty}(\Rd))}
			\lesssim \e .
		\end{split}
	\end{equation*}
	Similarly,
	$$
	\|\eta_2\|_{L^{1,\infty}(\Omega_T,\mu)}
	\lesssim \e .
	$$

	On the other hand, for every $\lambda>0$,
	\begin{equation*}
		\begin{split}
			&\mu\left(
			\left\{
			(t,x,y)\in\Omega_T:\
			|\eta_1(t,x)|
			>\lambda
			\right\}
			\right) \\
			&\leq
			\mu\left(
			\left\{
			(t,x,y)\in\Omega_T:\
			\frac{|\hat{b}(t,x)|}{\delta}>\lambda
			\right\}
			\right) \\
			&=
			\int_{\left\{(t,x,y)\in\Omega_T:\
				\frac{|\hat{b}(t,x)|}{\delta}>\lambda\right\}}
			\pi_t(\d x,\d y)\d t \\
			&\leq
			\int_{\left\{(t,x)\in[0,T]\times\Rd:\
				\frac{|\hat{b}(t,x)|}{\delta}>\lambda\right\}}
			\rho(t,x)\,\d x\d t \\
			&\leq
			\|\rho\|_{L^\infty([0,T]\times\Rd)}
			(\mathcal{L}^1\otimes\mathcal{L}^d)
			\left(
			\left\{
			(t,x)\in[0,T]\times\Rd:\
			\frac{|\hat{b}(t,x)|}{\delta}>\lambda
			\right\}
			\right).
		\end{split}
	\end{equation*}
	Hence
	$$
	\|\eta_1\|_{L^{p,\infty}(\Omega_T,\mu)}
	\lesssim
	\frac{1}{\delta}
	\|\hat{b}\|_{L^{p,\infty}([0,T]\times\Rd;\mathcal{L}^1\otimes\mathcal{L}^d)}
	\lesssim
	\frac{1}{\delta},
	$$
	and similarly,
	$$
	\|\eta_2\|_{L^{p,\infty}(\Omega_T,\mu)}
	\lesssim
	\frac{1}{\delta}.
	$$
	Combining the estimates on $\|\eta_i\|_{L^{1,\infty}(\Omega_T,\mu)}$ and $\|\eta_i\|_{L^{p,\infty}(\Omega_T,\mu)}$, $i=1,2$,
	by the interpolation inequality \eqref{interp 2}, we obtain
	\begin{align}\label{pf-J1-2}
		\|\eta_1\|_{L^1(\Omega_T,\mu)}
		+
		\|\eta_2\|_{L^1(\Omega_T,\mu)}
		\lesssim
		\e\left[1+\log\left(\frac{1}{\e\delta}\right)\right].
	\end{align}
	
	Combining the preceding estimates \eqref{pf-J1-1}, \eqref{pf-eta3} and \eqref{pf-J1-2} together, we derive
	\begin{align}\label{pf-J1}
		\mathrm{J}_1
		\lesssim
		\frac{1}{\delta}
		\int_0^T\int_{\Rd}
		|b(t,x)-\hat{b}(t,x)|\rho(t,x)\,\d x\d t
		+
		\e\left[1+\log\left(\frac{1}{\e\delta}\right)\right]
		+
		C_\e .
	\end{align}
	
	By \eqref{pf-se2-1} and the elementary inequality $\log(1+\sqrt{a^2+\delta^2}/\delta)\leq \log2+\log(1+a/\delta)$ for any $a\geq0$, we have
	\begin{align}\label{pf-se2-2}
		\E \log \left(1+\frac{\sqrt{|X_0-\hat{X}_0|^2+\delta^2}}{\delta}\right)
		\leq\mathcal{D}_{\delta}(\nu,\hat{\nu})+\log 2.
	\end{align}
	
	Therefore, the desired estimate follows by combining \eqref{pf-se2-2}, \eqref{pf-1} with  the estimates \eqref{pf-J1}  on $\mathrm{J}_1$ and
	\eqref{pf-J2} on $\mathrm{J}_2$, since $[0,\infty)\ni\alpha\mapsto \E \log [1+(|X_t-\hat{X}_t|^2+\alpha)^{1/2}/\delta]$ is increasing.
\end{proof}

\subsection{An application to well-posedness}\label{app}\hskip\parindent
In this subsection, we prove well-posedness of weak solutions to the FPE. Since the
condition on $\nabla b$ is insufficient for existence, we impose an additional
structural assumption on $b$ itself.

\begin{assumption}\label{assump 3}
	For each $j=1,\dots,d$, let $\Phi_j:\mathbb{R}^d\to\mathbb{R}$ be a homogeneous function of degree $1-d$, i.e.,
	\[
	\Phi_j(\lambda x)=\lambda^{1-d}\Phi_j(x),\qquad \lambda>0,\ x\in\mathbb{R}^d,
	\]
	which satisfies the following conditions:
	\begin{itemize}
		\item[]$(\Phi_1)$ $\Phi_j|_{\mathbb{R}^d\setminus\{0\}}\in C^2(\mathbb{R}^d\setminus\{0\})$;
		
		\item[]$(\Phi_2)$ there exists a constant $C>0$ such that
		for every integer $0\le m\le 2$,
		$$|\nabla^m\Phi_j(x)| \le\frac{ C }{|x|^{d-1+m}},\qquad x\neq 0,$$
		where $|\cdot|$ is understood as the Euclidean norm if $m=0,1$, and the Hilbert--Schmidt norm if $m=2$;

	\end{itemize}
	The drift vector field $b=(b_1,\dots,b_d)$ is defined componentwise by
	$$b_j(t,x):=\{\Phi_j * w_j(t,\cdot)\}(x),\qquad t\in[0,T],\, x\in\mathbb{R}^d,$$
	where $w=(w_1,\dots,w_d)\in L^\infty([0,T]; L^1(\mathbb{R}^d;\mathbb{R}^d))$.
\end{assumption}

\begin{remark}
	Assumption~\ref{assump 3} ensures that  for each $i,j=1,\cdots,d$,  $\partial_j\Phi_i$ is a singular kernel of fundamental type in the sense of Definition \ref{def:singular-kernel-fundamental-type}. Indeed, all the requirements in Definition~\ref{def:singular-kernel-fundamental-type} other than the cancellation
	condition follow immediately from the homogeneity and $(\Phi_1)$--$(\Phi_2)$ in Assumption \ref{assump 3}. It therefore remains only to verify the cancellation
	condition. Fix $i,j\in\{1,\ldots,d\}$, and let $e_j\in\mathbb{R}^d$ denote the $j$-th standard basis vector. For every $0<R_1<R_2<\infty$, since $\Phi_i|_{\mathbb{R}^d\setminus\{0\}}\in C^1(\mathbb{R}^d\setminus\{0\})$, Green's formula applied to the vector field $F(x):=\Phi_i(x)e_j$ on the annulus $B_{R_2}\setminus B_{R_1}$ gives
	\begin{align*}
		\int_{B_{R_2}\setminus B_{R_1}}
		\partial_j\Phi_i(x)\,\d x
		&=
		\int_{|x|=R_2}
		\Phi_i(x)\frac{x_j}{R_2}\,\d\omega(x)
		-
		\int_{|x|=R_1}
		\Phi_i(x)\frac{x_j}{R_1}\,\d\omega(x),
	\end{align*}
	where $\omega$ is the volume measure on the $(d-1)$-dimensional sphere. For every $0<R<\infty$, changing variables $x=Rz$ with $|z|=1$ and using the homogeneity of $\Phi_i$ of degree $1-d$, we have
	$$
	\int_{|x|=R}
	\Phi_i(x)\frac{x_j}{R}\,\d\omega(x)
	=
	\int_{|z|=1}
	\Phi_i(z)z_j\,\d\omega(z).
	$$
	Thus, the boundary integral is independent of $R$, and consequently
	$$
	\int_{B_{R_2}\setminus B_{R_1}}
	\partial_j\Phi_i(x)\,\d x=0,
	$$
	for every $0<R_1<R_2<\infty$.
\end{remark}

The main result in this subsection is contained in the following theorem.
\begin{theorem}\label{thm:well-posedness-fpe}
	Assume that Assumption \ref{assump 3}, $(\widetilde{\mathcal{H}}^{\sigma})$, $(\mathcal{H}^{b^\sigma})$, and $(\mathcal{H}^{\nu})$ hold. Then the FPE \eqref{FPE} admits a unique
	weak solution in the class $L^\infty([0,T]; L^1(\R^d)\cap L^\infty(\R^d))$.
\end{theorem}

Before proving the theorem, we present the following  basic regularity.
\begin{lemma}\label{regularity of b}
	Let $h  \in L^\infty([0,T]; L^1(\mathbb{R}^d))$ and let $\psi:\R^d\rightarrow\R$ satisfy
	$$|\psi(x)| \le\frac{ C }{|x|^{d-1}},\qquad x\neq 0,$$
	for some constant $C>0$. For $(t,x)\in[0,T]\times\R^d$, set $V(t,x):=[h(t,\cdot)*\psi](x)$.  Then
	\begin{equation}\label{b weak lebesgue}
		V\in L^{p,\infty}([0,T]\times \Rd),
		\qquad p:=\frac{d}{d-1}.
	\end{equation}
\end{lemma}
\begin{proof}
	By the assumption on $\psi$, we have the pointwise estimate
	$$|V(t,x)|\le \int_{\mathbb R^d} \frac{C |h(t,y)|}{|x-y|^{d-1}}\,\d y,\quad t\in[0,T],$$
	where $C>0$ is a constant. For any $R>0$ and a.e. $t\in[0,T]$,
	\begin{align*}
		|V(t,x)|&\le C\int_{B_R(x)} \frac{|h(t,y)|}{|x-y|^{d-1}}\,\d y
		+ C\int_{\mathbb R^d\setminus B_R(x)} \frac{|h(t,y)|}{|x-y|^{d-1}}\,\d y \\
		&\le C\left(\mathbbm{1}_{B_R}\frac{1}{|\cdot|^{d-1}}\right)*|h(t,\cdot)|(x)
		+ \frac{C N}{R^{d-1}},
	\end{align*}
	for some constant $N>0$, where the last inequality is due to that $\|h(t,\cdot)\|_{L^1}\le N$ for a.e. $t\in[0,T]$ since $h\in L^\infty([0,T];L^1(\mathbb R^d))$.
	
	Let $\lambda>0$ be arbitrary and choose $R=(2CN/\lambda)^{1/(d-1)}$, so that $CN/R^{d-1}=\lambda/2$. Letting $\mu=\mathcal L^1|_{[0,T]}\otimes\mathcal L^d$, by Chebyshev's inequality, we derive
	\begin{align*}
		\mu(\{|V|>\lambda\})
		&\le \mu\left(\left\{ \left(\mathbbm{1}_{B_R}\frac{1}{|\cdot|^{d-1}}\right)*|h| > \frac{\lambda}{2C} \right\}\right) \\
		&\le \frac{2C}{\lambda}\int_0^T\int_{\mathbb R^d}\left(\mathbbm{1}_{B_R}\frac{1}{|\cdot|^{d-1}}\right)*|h(t,\cdot)|(x)\,\d x \d t.
	\end{align*}
	By Young's convolution inequality and the computation
	$$\left\|\mathbbm{1}_{B_R}\frac{1}{|\cdot|^{d-1}}\right\|_{L^1(\R^d)}= \mathfrak{c}_d R,$$
	we obtain
	$$\mu(\{|V|>\lambda\})\le \frac{2C}{\lambda} \cdot (\mathfrak{c}_d R) \cdot (TN)\le C' \lambda^{-\frac{d}{d-1}},$$
	for some constant $C'>0$, where $\mathfrak{c}_d$ is the volume of the unit ball in $\R^d$.
	Therefore, with $p=\frac{d}{d-1}$,
	$$\sup_{\lambda>0}\lambda^p\mu(\{|V|>\lambda\})<\infty,$$
	which proves \eqref{b weak lebesgue}.
\end{proof}

The next lemma computes the distributional derivative of $b$.
\begin{lemma}[Distributional derivative of $b$]
	Let Assumption \ref{assump 3} hold. Then, for every $t\in[0,T]$, the
	distributional derivative of $b(t,\cdot)$ is given by
	\begin{align}\label{derivative-b}
		\partial_j b_i(t,x)
		=[\partial_j\Phi_i*w_i(t,\cdot)](x)-c_{ij}w_i(t,x)
		\qquad \text{in } \mathscr{D}'(\Rd),
	\end{align}
	where the convolution is understood in the distributional sense of \eqref{distribution extension} and $c=(c_{ij})_{1\leq i,j\leq d}$ is the constant matrix with
	$$	c_{ij}
	=
	\int_{|x|=1}\Phi_i(x)x_j\,\d \omega(x),
	\qquad 1\leq i,j\leq d.
	$$
\end{lemma}
\begin{proof}
	Fix $i,j\in\{1,\cdots,d\}$ and take an arbitrary $\varphi\in C_c^\infty(\Rd)$. By assumptions $(\Phi_1)$ and $(\Phi_2)$ with $m=0$, we have $\Phi_i\in L^1_{\rm loc}(\mathbb R^d)$. For every $\epsilon>0$, applying Green's formula on $\{|x|\ge\epsilon\}$ yields
	\begin{align}\label{pf-derivative-b-1}
		\int_{|x|\ge\epsilon}\Phi_i(x)\partial_j\varphi(x)\,\d x=
		-\int_{|x|\ge\epsilon}\varphi(x)\partial_j\Phi_i(x)\,\d x  + \int_{|x|=\epsilon}\Phi_i(x)\varphi(x)\frac{x_j}{|x|}\,\d\omega(x).
	\end{align}
	Since $\Phi_i\partial_j\varphi\in L^1(\R^d)$, it is clear that letting $\epsilon\to0^+$, the left-hand side of \eqref{pf-derivative-b-1} converges to $\int_{\mathbb R^d}\Phi_i(x)\partial_j\varphi(x)\,\d x$. For the boundary term, changing variables $x=\epsilon z$ with $|z|=1$ and using the homogeneity $\Phi_i(\epsilon z)=\epsilon^{1-d}\Phi_i(z)$, we obtain
	\begin{align*}
		\int_{|x|=\epsilon}\Phi_i(x)\varphi(x)\frac{x_j}{|x|}\,d\omega(x)
		&=\int_{|z|=1}\Phi_i(z)\varphi(\epsilon z)z_j\,d\omega(z)\\
		&\rightarrow\varphi(0)\int_{|z|=1}\Phi_i(z)z_j\,d\omega(z),\quad\mbox{as }\epsilon\rightarrow0^+,
	\end{align*}
	where convergence follows from the dominated convergence theorem, since $\varphi(\epsilon z)\to\varphi(0)$ uniformly and $\Phi_i(z)z_j$ is integrable on $\{|z|=1\}$. Consequently, the limit of the first integral on the right-hand side of \eqref{pf-derivative-b-1} exists, and it is precisely the Cauchy principal value
	$$\mathrm{p.v.}\int_{\mathbb R^d}\varphi(x)\partial_j\Phi_i(x)\,\d x:=\lim_{\epsilon\to0^+}\int_{|x|\ge\epsilon}\varphi(x)\partial_j\Phi_i(x)\,\d x.$$
	Thus,
	\begin{align*}
		\int_{\mathbb R^d}\Phi_i(x)\partial_j\varphi(x)\,\d x=
		-\mathrm{p.v.}\int_{\mathbb R^d}\varphi(x)\partial_j\Phi_i(x)\,\d x + c_{ij}\phi(0),
	\end{align*}
	which in the sense of distributions reads
	\begin{align}\label{pf-derivative-b-2}
		\partial_j\Phi_i = {\rm p.v.}(\partial_j\Phi_i) - c_{ij}\delta_0.
	\end{align}
	where  $\delta_0$ is the Dirac mass at the origin, and
	$$c_{ij}:=\int_{|x|=1}\Phi_i(x)x_j\,\d \omega(x).$$

	Fix $t\in[0,T]$. Recall that
	\begin{align*}
		\<\partial_j\Phi_i*w_i(t,\cdot),\varphi\>
		:=
		\<w_i(t,\cdot),\widetilde{\partial_j\Phi_i}*\varphi\>,
	\end{align*}
	where
	$$
	\widetilde{\partial_j\Phi_i}(x)
	:=
	\partial_j\Phi_i(-x),
	\qquad x\in\Rd .
	$$
	Using \eqref{pf-derivative-b-2}, we get
	\begin{align*}
		\<b_i(t,\cdot),\partial_j\varphi\>
		&=
		\int_{\Rd}\int_{\Rd}
		\Phi_i(x-y)w_i(t,y)\partial_j\varphi(x)\,\d y\d x \\
		&=\int_{\Rd}\left(\int_{\Rd}
		\Phi_i(x-y)\partial_j\varphi(x)\,\d x\right)w_i(t,y)\,\d y\\
		&=
		\int_{\Rd}
		\left(
		-\operatorname{p.v.}\int_{\Rd}
		\varphi(x)\partial_j\Phi_i(x-y)\,\d x
		+
		c_{ij}\<\delta_0(\cdot-y),\varphi\>
		\right)
		w_i(t,y)\,\d y \\
		&=
		-\<\partial_j\Phi_i*w_i(t,\cdot),\varphi\>+\<c_{ij}w_i(t,\cdot),\varphi\>.
	\end{align*}
	Therefore
	$$
	\partial_j b_i(t,x)
	=
	\{\partial_j\Phi_i*w_i(t,\cdot)\}(x) -
	c_{ij}w_i(t,x)
	\qquad \text{in } \mathscr{D}'(\Rd),
	$$
	which proves the lemma.
\end{proof}

The zeroth-order term on the right-hand side of \eqref{derivative-b}, which originates from the Dirac mass, is more regular than the singular integral, and it is fully controlled by the estimates for the singular part. Hence, for notational convenience and without loss of generality, we write
$$
\nabla b(t,x)=\{\nabla\Phi*w(t,\cdot)\}(x),
\qquad \text{in } \mathscr{D}'(\Rd;\R^{d\times d}),
$$
for any $t\in[0,T]$, where the convolution is understood in the sense of $(\widetilde{\mathcal{H}}^{\nabla b})$.

\begin{proof}[Proof of Theorem \ref{thm:well-posedness-fpe}]
	By Corollary \ref{cor:uniqueness}, weak solutions in the class
	$L_t^\infty(L_x^1\cap L_x^\infty)$ are unique. It therefore remains only to
	prove existence.
	
	Let $(\phi^\e)_{\e\in(0,1)}$ be a standard mollifier defined by \eqref{mollifier}. Set
	$$
	b^\e:=b*\phi^\e,
	\quad
	a^\e:=a*\phi^\e\in L^1\big([0,T];C_b^\infty(\Rd;\R^{d\times d})\big),
	\quad
	\rho_0^\e:=\rho_0*\phi^\e\in C_b^\infty(\Rd).
	$$
	By Lemma \ref{regularity of b} and the interpolation inequality \eqref{interp 1},
	$$
	b\in L_{\mathrm{loc}}^q([0,T]\times\Rd;\R^d),
	\qquad
	1\leq q<\frac{d}{d-1}.
	$$
	Hence $b^\e(t,\cdot)\in C^\infty(\Rd;\R^d)$ for every 
	$\e\in(0,1)$ and $t\in(0,T)$.
	Moreover, Assumption
	\ref{assump 3} gives
	$$
	b^\e\in L^1([0,T];C_b^\infty(\Rd;\R^d)).
	$$
	
	Repeating the proof of Theorem \ref{thm partial}(1), it remains
	only to verify that
	$$
	\sup_{0<\e<1}\int_0^T\int_{\Rd}
	|b^\e(s,x)|\rho^\e(s,x)\,\d x\,\d s<\infty.
	$$
	Set
	$$
	\mu^\e(\d t,\d x)
	:=
	\mathbbm{1}_{[0,T]}(t)\rho^\e(t,x)\,\d x\d t,\qquad \mu:=\mathcal{L}^1|_{[0,T]}\otimes\mathcal{L}^d.
	$$
	By \eqref{interp 1} and the uniform boundedness of $\rho^\e$, we obtain, for
	$p=\frac{d}{d-1}$,
	\begin{equation*}
		\begin{split}
			\|b^\e\|_{L^1([0,T]\times \Rd;\mu^\e)}
			&\leq
			\frac{p}{p-1}T^{1-1/p}
			\|b^\e\|_{L^{p,\infty}([0,T]\times\Rd;\mu^\e)}\\
			&\lesssim
			\|b^\e\|_{L^{p,\infty}([0,T]\times \Rd;
				\mu)}.
		\end{split}
	\end{equation*}
	It remains to control the weak $L^{p}$ quasi-norm of $b^\e$. We claim that
	$$
	\|b^\e\|_{L^{p,\infty}([0,T]\times \Rd;
		\mu)}
	\lesssim
	\|b\|_{L^{p,\infty}([0,T]\times \Rd;
		\mu)},\quad\e\in(0,1).
	$$
	Indeed, for $\lambda>0$, write
	$$
	b=b_1+b_2,
	\qquad
	b_1:=b\,\mathbbm{1}_{\{|b|>\lambda/2\}},
	\qquad
	b_2:=b\,\mathbbm{1}_{\{|b|\leq \lambda/2\}}.
	$$
	Let $b_1^\e=b_1*\phi^\epsilon$ and $b_2^\e=b_2*\phi^\epsilon$. Then $|b_2^\e|\leq \lambda/2$, and hence
	$$
	\{|b^\e|>\lambda\}\subset \{|b_1^\e|>\lambda/2\}.
	$$
	Thus,
	$$
	\mu(\{|b^\e|>\lambda\})
	\leq
	\frac{2}{\lambda}\|b_1^\e\|_{L^1_{t,x}}
	\leq
	\frac{2}{\lambda}\|b_1\|_{L^1_{t,x}}.
	$$
	Since
	$b\in L^{p,\infty}([0,T]\times \Rd; \R^d)$ with $p=\frac{d}{d-1}>1$ by Lemma \ref{regularity of b}, applying the layer-cake formula, we obtain
	\begin{align*}
		\|b_1\|_{L^1_{t,x}}
		&=
		\int_0^\infty
		\mu\left(\{|b_1|>s\}\right)\,\d s\\
		&=\int_0^{\frac{\lambda}{2}}
		\mu(\{|b_1|>s\})\,\d s +  \int_{\frac{\lambda}{2}}^\infty 
		\mu(\{|b_1|>s\})\,\d s\\
		&=
		\frac{\lambda}{2}
		\mu\left(\left\{|b|>\frac{\lambda}{2}\right\}\right)
		+
		\int_{\frac{\lambda}{2}}^\infty
		\mu(\{|b|>s\})\,\d s \\
		&\leq
		\frac{\lambda}{2}
		\|b\|_{L^{p,\infty}([0,T]\times \Rd;
			\mu)}^p
		\left(\frac{\lambda}{2}\right)^{-p}
		+
		\|b\|_{L^{p,\infty}([0,T]\times \Rd;
			\mu)}^p
		\int_{\frac{\lambda}{2}}^\infty s^{-p}\,\d s \\
		&=
		\|b\|_{L^{p,\infty}([0,T]\times \Rd;
			\mu)}^p
		\left(\frac{\lambda}{2}\right)^{1-p}
		+
		\frac{1}{p-1}
		\|b\|_{L^{p,\infty}([0,T]\times \Rd;
			\mu)}^p
		\left(\frac{\lambda}{2}\right)^{1-p} \\
		&=
		\frac{p\,2^{p-1}}{p-1}
		\|b\|_{L^{p,\infty}([0,T]\times \Rd;
			\mu)}^p
		\lambda^{1-p}.
	\end{align*}
	It follows that
	$$
	\mu(\{|b^\e|>\lambda\})
	\leq
	\frac{p\,2^p}{p-1}
	\|b\|_{L^{p,\infty}([0,T]\times \Rd;
		\mu)}^p
	\lambda^{-p}.
	$$
	Equivalently,
	$$
	\lambda
	\bigl[
	\mu(\{|b^\e|>\lambda\})
	\bigr]^{1/p}
	\leq
	2\left(\frac{p}{p-1}\right)^{1/p}
	\|b\|_{L^{p,\infty}([0,T]\times \Rd;
		\mu)} .
	$$
	Taking the supremum over $\lambda>0$ gives
	$$
	\|b^\e\|_{L^{p,\infty}([0,T]\times\Rd;
		\mu)}
	\leq
	2\left(\frac{p}{p-1}\right)^{1/p}
	\|b\|_{L^{p,\infty}([0,T]\times\Rd;
		\mu)}.
	$$
	
	Arguing as in the 
	proof of Theorem \ref{thm partial}(1),
	we obtain
	$$
	\rho\in L_t^\infty(L_x^1\cap L_x^\infty),
	$$
	and $\rho(t,\cdot)$ is a probability density for a.e. $t\in[0,T]$.
	Moreover, $\rho$ satisfies \eqref{weak solution integral}. Since the drift $b$
	belongs only to a weak Lebesgue space, it remains to verify the integrability
	condition \eqref{integral condition}. As $\sigma$ is uniformly bounded, it
	suffices to show that
	$$
	\|b\rho\|_{L^1([0,T]\times \Rd;\R^d)}<\infty.
	$$
	This follows from \eqref{b weak lebesgue} and the interpolation inequality
	\eqref{interp 1}.
	
	The proof is complete.
\end{proof}

			\section*{Acknowledgment}
			The authors gratefully acknowledge the financial support of the National Key R\&D Program of China (No. 2022YFA1006000).
			\section*{Data Availability Statement}
			 Data sharing is not applicable to this article as no datasets were generated or
			analysed during the current study.
			\section*{Declarations}
			
			\noindent\textbf{Conflict of Interests.}
			The authors declare that they have no conflict of interest.
			
			\medskip
			
			\noindent\textbf{Declaration of generative AI and AI-assisted technologies in the manuscript preparation process.}
			During the preparation of this work, the authors used ChatGPT (OpenAI) solely for language editing and polishing. It was not used for research design, mathematical analysis, proofs, interpretation of results, or conclusions. The authors reviewed and edited all relevant content and take full responsibility for the published article.


\begin{thebibliography}{4}
				\bibitem{Ambrosio 2004}
			L. Ambrosio, Transport equation and Cauchy problem for BV vector fields, Invent. Math. 158 (2004),  227--260.
			
			\bibitem{Ambrosio 2008}
			L. Ambrosio, N. Gigli, G. Savar\'e, \emph{Gradient flows in metric spaces and in
				the space of probability measures}, 2nd
			ed., Lectures in Mathematics ETH Z\"urich, Birkh\"auser Verlag, Basel, 2008.
			
			\bibitem{Bogachev 2022}
			V. I. Bogachev, N. V. Krylov, M. R\"{o}ckner, S. V. Shaposhnikov, \textit{Fokker--Planck--Kolmogorov Equations,} Vol. 207, American Mathematical Society, Providence, RI, 2022.
			
			\bibitem{Bouchut 2016}
			A. Bohun, F. Bouchut,  G. Crippa, Lagrangian flows for vector fields with anisotropic regularity, Ann. Inst. H. Poincar\'e C Anal. Non Lin\'eaire. 33 (6) (2016), 1409--1429.
			
			\bibitem{Brezis 2011}
			H. Br\'ezis, \textit{Functional Analysis, Sobolev Spaces and Partial
				Differential Equations,} Springer, New York, 2011.
			
			\bibitem{Bouchut 2013}
			F. Bouchut, G. Crippa, Lagrangian flows for vector fields with gradient given by a singular integral, J. Hyperbolic Differ. Equ. 10 (2) (2013),  235--282.
			
			\bibitem{Crippa 2020}
			G. Crippa, S. Ligabue, A note on the Lagrangian flow associated to a partially regular vector field, Differ. Equ. Dyn. Syst. 31 (4) (2020), 767--786.
			
			\bibitem{Crippa 2017}
			G. Crippa, C. Nobili, C. Seis, and S. Spirito, Eulerian and Lagrangian solutions to the continuity and Euler equations with $L^1$-vorticity, SIAM J. Math. Anal. 49 (5) (2017), 3973--3998.
			
			\bibitem{Crippa 2008}
			G. Crippa, C. De Lellis, Estimates and regularity results for the DiPerna--Lions flow, J. Reine Angew. Math. 616 (2008), 15--46.
			
			\bibitem{DiPerna 1989}
			R. J. DiPerna, P.-L. Lions, Ordinary differential equations, transport theory and Sobolev spaces, Invent. Math. 98 (1989), no. 3, 511--547.
			
			\bibitem{DiPernaLions1988 boltzmann}
			R. J. DiPerna, P.-L. Lions, On the Fokker--Planck--Boltzmann equation,
			Comm. Math. Phys. 120 (1988), 1--23.
			
			\bibitem{DiPernaLions1989 boltzmann}
			R. J. DiPerna, P.-L. Lions, On the Cauchy problem for Boltzmann equations:
			global existence and weak stability, Ann. of Math. 130 (1989), no. 2,
			321--366.
			
			\bibitem{Evans 2025}
			L. C. Evans, \textit{Measure theory and fine properties of functions}, 2nd ed., Chapman and
			Hall/CRC, New York, 2025.
			
			\bibitem{Fang 2010}
			S. Fang, D. Luo, A. Thalmaier, Stochastic differential equations with coefficients in Sobolev spaces, J. Funct. Anal. 259 (5) (2010), 1129--1168.
			
			\bibitem{Flandoli 2010}
			F. Flandoli, M. Gubinelli, E. Priola, Well-posedness of the transport equation by stochastic perturbation. Invent. Math. 180 (2010), 1--53.
			
			\bibitem{Figalli 2008}
			A. Figalli, Existence and uniqueness of martingale solutions for SDEs with
			rough or degenerate coefficients, J. Funct. Anal. 254 (2008), 109--153.
			
			
			\bibitem{Grafakos2014a}
			L. Grafakos, \emph{Classical Fourier analysis}. 3rd ed., Graduate Texts in Mathematics, 249. Springer, New York, 2014.
			
			\bibitem{Ikeda 1989}
			N. Ikeda, S. Watanabe, \textit{Stochastic Differential Equations and Diffusion Processes}, 2nd
			ed., North-Holland Mathematical Library, vol. 24, North-Holland, Amsterdam, 1989.
			
			\bibitem{Jabin2010}
			P. E. Jabin, Differential equations with singular fields, J. Math. Pures Appl. (9) 94 (2010), 597--621.
			
			\bibitem{LeBLions19}
			C. Le Bris, P.-L. Lions, \emph{Parabolic equations with irregular data and related issues: applications to stochastic differential equations}, De Gruyter Series in
			Applied and Numerical Mathematics, vol 4, De Gruyter, Berlin, 2019.
			
			\bibitem{Li 2019}
			H. Li, D. Luo, Quantitative stability estimates for Fokker--Planck equations, J. Math. Pures Appl. (9) 122  (2019), 125--163.
			
			\bibitem{Luo2013}	
			D. Luo, Fokker--Planck type equations with Sobolev diffusion coefficients and BV drift coefficients, Acta Math. Sin.
			Engl. Ser. (2) 29 (2013), 303--314.
			
			\bibitem{Luo 2014}
			D. Luo, Uniqueness of degenerate Fokker--Planck equations with weakly differentiable drift whose gradient is given by a singular integral.
			Electron. Commun. Probab. 19  (2014), no. 43, 1--14.
			
			\bibitem{Lions 2008}
			C. Le Bris, P.L. Lions, Existence and uniqueness of solutions to Fokker--Planck type equations with irregular coefficients, Commun. Partial Differ. Equ. 33 (2008), 1272--1317.
			
			\bibitem{Le Bris 2004}
			C. Le Bris, P.L. Lions, Renormalized solutions of some transport equations with partially $W^{1,1}$ velocities and applications, Ann. Mat. Pura Appl. 183 (2004), 97--130.
			
			\bibitem{Ladyzenskaja 1968}
			O. A. Lady{\v{z}}enskaja, V. A. Solonnikov, 
			N. N. Ural'ceva, \textit{Linear and quasilinear equations of parabolic type}, Translations of Mathematical Monographs, vol. 23,
			American Mathematical Society, Providence, R.I., 1968.
			
			
			
			\bibitem{Nguyen 2021}
			Q.-H. Nguyen, Quantitative estimates for regular Lagrangian flows with BV vector fields, Comm. Pure Appl. Math. 74 (2021), no. 6, 1129--1192.
			
			\bibitem{Seis 2022}
			V. Navarro-Fern\'{a}ndez, A. Schlichting, C. Seis, Optimal stability estimates and a new uniqueness result for advection--diffusion equations, Pure Appl. Anal. 4 (3) (2022), 571--596.
			
			
			\bibitem{Risken 1989}
			H. Risken and T. Frank, \textit{The Fokker--Planck Equation: Methods of Solution and Applications,} 2nd ed., Springer, Berlin, 1989.
			
			
			
			\bibitem{Rockner 2010}
			M. R\"ockner, X. Zhang, Weak uniqueness of Fokker--Planck equations with degenerate and bounded coefficients,
			C. R. Math. Acad. Sci. Paris. 348, (2010) 435--438.
			
			\bibitem{StroockVaradhan 2006}
			D. W. Stroock and S. R. S. Varadhan, \textit{Multidimensional Diffusion
				Processes,} Classics in Mathematics, Springer, Berlin, 2006. Reprint of
			the 1997 edition. 
			
			\bibitem{Seis 2017}
			C. Seis, A quantitative theory for the continuity equation, Ann. Inst. Henri Poincar\'e, Anal. Non Lin\'eaire. 34, (2017)
			1837--1850.
			
			\bibitem{Trevisan 2016}
			D. Trevisan, Well-posedness of multidimensional diffusion processes with weakly differentiable coefficients, Electron. J. Probab. (22) 21, (2016), 41 pp.
				
			\end{thebibliography}
			\end{document}